\documentclass[twocolumn]{autart}

\usepackage[utf8]{inputenc}
\usepackage[T1]{fontenc}

\usepackage{amsmath,amsfonts,amssymb} 
\usepackage{mathtools}
\usepackage{mathabx}
\usepackage{mathrsfs}
\usepackage{dsfont}
\usepackage{scalerel}
\usepackage{listings}
\usepackage{graphicx}
\usepackage{subcaption}
\usepackage{enumerate}
\usepackage{cite}
\usepackage{url}
\usepackage{csquotes}
\usepackage{xcolor}
\usepackage{multicol}
\usepackage{multirow}
\usepackage{makecell}
\usepackage{soul}

\usepackage{graphicx}
\graphicspath{ {Figures/} }
\DeclareGraphicsExtensions{.eps}
\usepackage{subcaption}
\usepackage{calc}
\usepackage{tikz,pgfplots}
\usetikzlibrary{shapes, decorations.markings, calc}
\tikzstyle{vertex}=[circle, draw, inner sep=0pt, minimum size=6pt]
\newcommand{\vertex}{\node[vertex]}
\usetikzlibrary{positioning,hobby,arrows,calc}
\tikzset{
modal/.style={>=stealth',shorten >=1pt,shorten <=1pt,auto,node distance=1.5cm,
semithick},
world/.style={circle,draw,minimum size=0.5cm,fill=gray!15},
point/.style={circle,draw,inner sep=0.5mm,fill=black},
reflexive above/.style={->,loop,looseness=7,in=100,out=60},
reflexive below/.style={->,loop,looseness=7,in=250,out=290},
reflexive left/.style={->,loop,looseness=7,in=150,out=200},
reflexive right/.style={->,loop,looseness=7,in=30,out=340}
}

\usepackage{algorithm}

\let\classAND\AND
\let\AND\relax
\usepackage{algorithmic}

\let\AND\classAND
\AtBeginEnvironment{algorithmic}{\let\AND\algoAND}

\floatname{algorithm}{Algorithm}

\newtheorem{Theorem}{Theorem}

\newtheorem{Corollary}{Corollary}[Theorem]
\newtheorem{Lemma}[Theorem]{Lemma}

\newtheorem{Assumption}{Assumption}

\definecolor{blue1}{RGB}{128, 191, 255}
\definecolor{blue2}{RGB}{51   153  255}
\definecolor{blue3}{RGB}{50,60,200}
\definecolor{blue4}{RGB}{30,120,200}
\definecolor{redd1}{RGB}{255, 153, 128}
\definecolor{redd2}{RGB}{255  71   26}
\definecolor{redd3}{RGB}{200,50,50}
\definecolor{gren1}{RGB}{153, 230, 153}
\definecolor{gren2}{RGB}{51   204  51}
\definecolor{gren3}{RGB}{50,140,50}
\definecolor{yelw1}{RGB}{255, 179, 26}
\definecolor{yelw2}{RGB}{255  204  0}
\definecolor{yelw3}{RGB}{241,188,49}
\definecolor{prpl1}{RGB}{204, 153, 255}
\definecolor{prpl2}{RGB}{153  51   255}
\definecolor{mage1}{RGB}{255, 51, 51}
\definecolor{gray1}{RGB}{150,150,150}
\definecolor{brwn1}{RGB}{120,70,20}

\def \ba{\begin{array}}
\def \ea{\end{array}}
\def \bea{\begin{eqnarray}}
\def \eea{\end{eqnarray}}
\def \be{\begin{equation}}
\def \ee{\end{equation}}
\def \BEA{\begin{eqnarray*}}
\def \EEA{\end{eqnarray*}}
\def \BE{\begin{equation*}}
\def \EE{\end{equation*}}
\def \disp{\displaystyle}
\def \colsep{\arraycolsep}
\def \mb{\mathbf}
\def \bb{\mathbb}
\def \bs{\boldsymbol}
\def \mc{\mathcal}
\def \mf{\mathfrak}

\def \mr{\mathring}
\def \diag{\mathtt{diag}}
\def \rank{\mathtt{rank}}
\def \grank{\mathtt{grank}}

\def \eig{\mathtt{eig}}

\def \trace{\mathtt{trace}}
\def \sym{\mathtt{sym}}
\def \Ker{\mathtt{ker}}

\def \-{\dagger} 

\def \T{\intercal}

\def \a{\mathtt{a}}

\def \ol{\overline}
\def \ul{\underline}

\usepackage{textcomp}
\def\BibTeX{{\rm B\Kern-.05em{\sc i\Kern-.025em b}\Kern-.08em
    T\Kern-.1667em\lower.7ex\hbox{E}\Kern-.125emX}}

\begin{document}

\begin{frontmatter}

\title{
Clustering-Based Average State Observer Design for Large-Scale Network Systems
}

\author[umar]{Muhammad Umar B. Niazi}\ead{mubniazi@kth.se},
\author[xiaodong]{Xiaodong Cheng}\ead{xc336@cam.ac.uk},
\author[carlos]{Carlos Canudas-de-Wit}\ead{carlos.canudas-de-wit@gipsa-lab.fr},
\author[jacquelien]{Jacquelien M. A. Scherpen}\ead{j.m.a.scherpen@rug.nl}

\address[umar]{Division of Decision and Control Systems, EECS, KTH Royal Institute of Technology, SE-100 44 Stockholm, Sweden}
\address[xiaodong]{Department of Engineering, University of Cambridge, Cambridge CB2 1PZ, U.K.}
\address[carlos]{GIPSA-lab, CNRS, 38000 Grenoble, France}
\address[jacquelien]{Jan C. Willems Center for Systems and Control, ENTEG, Faculty of Science and Engineering, University of Groningen, 9747 AG Groningen, The Netherlands}

\thanks{This work was supported by the European Research Council through the European Union's Horizon 2020 Research and Innovation Programme (Scale-FreeBack) under Grant ERC-AdG no. 694209.}


\begin{abstract}
This paper addresses the aggregated monitoring problem for large-scale network systems with a few dedicated sensors. Full state estimation of such systems is often infeasible due to unobservability and/or computational infeasibility. Therefore, through clustering and aggregation, a tractable representation of a network system, called a projected network system, is obtained for designing a minimum-order average state observer. This observer estimates the average states of the clusters, which are identified with explicit consideration to the estimation error. Moreover, given the clustering, the proposed observer design algorithm exploits the structure of the estimation error dynamics to achieve computational tractability. Simulations show that the computation of the proposed algorithm is significantly faster than the usual $\mc{H}_2/\mc{H}_\infty$ observer design techniques. On the other hand, compromise on the estimation error characteristics is shown to be marginal.
\end{abstract}


\end{frontmatter}

\section{Introduction}

Knowledge of the system's state is undoubtedly crucial for monitoring and control. However, for large-scale network systems, it is challenging to estimate the full state vector with limited computational and sensing resources \cite{antoulas2005}. This is because low computational power results in the intractability of the state estimation algorithm, and few sensors render the network system unobservable.

With a limited number of sensors, it is reasonable to perform aggregated monitoring by clustering a large-scale network system and estimating the average states of the clusters. Real-world applications include urban traffic networks, building thermal systems, and water distribution networks. For instance, estimating average traffic densities in multiple sectors of an urban traffic network allows for monitoring congestion in different areas of the city~\cite{rodriguez-vega2021trb}, or estimating mean operative temperatures of rooms in a building allows for monitoring thermal comfort in its interior~\cite{niazi2020ifac, deng2010}. With such information on the state, it is then possible to control the network system in an aggregated sense \cite{nikitin2021tac}.

\subsection{Literature review}

Estimating the average states is the same as estimating linear functionals of the state, which has been a topic of interest for several decades. Luenberger \cite{luenberger1971} was the first to propose a linear functional observer with order equal to the system's observability index multiplied by the number of functionals to be estimated. Later, \cite{murdoch1973} and \cite{aldeen1999} showed that Luenberger's functional observer is conservative and its order can be reduced significantly. Approaching the problem from a different perspective, Darouach \cite{darouach2000} provided a necessary and sufficient condition for the existence and design of a functional observer with order equal to the number of linear functionals to be estimated, which is the minimum achievable order. The design procedure of Darouach's functional observer is recently improved in \cite{darouach2019}; however, the existence conditions are still restrictive, and the convergence is not always guaranteed, as the order of the functional observer is bounded by the number of functionals. This led to the development of the notion of functional observability in \cite{fernando2010} and \cite{jennings2011} followed by \cite{rotella2011, rotella2015, rotella2016}, which propose different methodologies to increase the order of Darouach's functional observer by a minimal amount to attain convergence. These methods are iterative and require rank computations of the concatenation of multiple observability matrices at every iteration, which becomes intractable for large-scale network systems. Therefore, for designing average state observers, it is requisite to cluster and aggregate the network system for obtaining a projected network system, which is an aggregated, tractable representation.

Methods based on the aggregation of large-scale systems also have a rich history. Having its roots in chemical reaction systems \cite{wei1969}, the notion of lumpability, which allows for an exact aggregated representation of a large-scale system, is studied rigorously in \cite{coxson1984} and \cite{atay2017}. For network systems, \cite{ji2007} showed that lumpability is equivalent to having an equitable partition of the underlying network. Later, for studying average controllability, the condition of equitable partition was relaxed to almost equitable partition in \cite{martini2010, egerstedt2012, aguilar2017}. However, under the constraints on sensor locations, the number of clusters, and cluster connectivity, achieving equitable or almost equitable partitions turn out to be very challenging~\cite{martin2019, niazi2019cdc}.

A similar line of research employs clustering or projection-based model reduction methods \cite{monshizadeh2014, ishizaki2014, ishizaki2015, cheng2017, cheng2019weakly} for approximating the average states. These model reduction tools not only preserve some dynamical properties of a network system but also its topological structure. Preserving the structure is important as monitoring and control of network systems usually rely on its underlying graph structure \cite{cheng2021}. The goal of clustering-based model reduction is to reduce the system's dimension by identifying and aggregating the clusters in a network system that yield minimum model approximation error, which is characterized in terms of $\mc{H}_2$ or $\mc{H}_\infty$ norms. In other words, the idea is to obtain a reduced, aggregated system with a tractable dimension whose input-output behavior is similar to the input-output behavior of the original network system. This allows for aggregated monitoring as the reduced system can be employed for estimating the approximated average states of the original network system. 

In this regard, \cite{sadamoto2017} presents an average state observer design based on the model reduction techniques developed in \cite{ishizaki2014, ishizaki2015}. A similar technique has been used to design an average Kalman filter in \cite{watanabe2015}. However, the design procedure is based on the solution of a Linear Matrix Inequality (LMI), which is not only computationally expensive but also doesn't provide an understanding on how inter-cluster and intra-cluster topologies affect the performance of an average state estimation algorithm. This motivated the development of the notions of average observability and average detectability in \cite{niazi2020tcns}, which provide the corresponding necessary and sufficient conditions on the inter-cluster and intra-cluster topologies of the network system. Since then, several average state observer designs have been proposed, for example, sliding mode design in \cite{pilloni2021} and design in the presence of outlier nodes in \cite{pratap2020, pratap2021}.

\subsection{Our contribution}

In this paper we present a clustering-based method to aggregate the network system and design an average state observer yielding a minimal asymptotic average state estimation error. The form of the average state observer is chosen to be similar to Darouach's functional observer \cite{darouach1994, darouach2000, darouach2019} with order equal to the number of clusters; however, the design criteria is adopted from \cite{niazi2020tcns}. 

The approach presented in this paper improves upon the work of \cite{sadamoto2017, niazi2020tcns}. Unlike \cite{niazi2020tcns}, we do not assume pre-specified clustering of a network system. We propose a cyclic coordinate descent scheme to achieve a suboptimal clustering-based average state observer. Such an approach is computationally efficient for large-scale network systems as compared to the LMI-based approach of \cite{sadamoto2017}. Moreover, the clustering method in \cite{sadamoto2017} doesn't consider an upper bound on the number of clusters, which may result in an infeasible solution. To address this issue, we consider a fixed number of clusters.

The solution to the clustering-based average state observer design problem naturally comprises two parts: finding an optimal clustering and finding an optimal average state observer. The proposed algorithm, therefore, iteratively seeks one by fixing the other. That is, given an initial clustering, we find the optimal average state observer design. Then, fixing the optimal average state observer design, we find the optimal clustering. This process is repeated until convergence or maximum number of iterations is reached. However, finding an optimal clustering is a non-convex, mixed integer-type optimization problem, which is an NP-hard problem \cite{burer2012}. Thus, a greedy clustering algorithm is proposed to obtain a suboptimal solution. On the other hand, finding an optimal average state observer is equivalent to $\mc{H}_2$ design, which is a convex optimization problem whose solution can be obtained through LMI formulation \cite{boyd2004, duan2013}. However, for large-scale network systems, we show that solving an LMI feasibility problem is computationally expensive. Thus, a structural relaxation on the average state observer design is required to achieve computational tractability. 

We provide a sufficient condition for the stabilizability of average state observer under structural relaxation and indicate its implications in the clustering, which are then integrated in the algorithm as clustering constraints.
Through a simulation example, we show that the computational time under our design is improved significantly, whereas the compromise on the optimality as compared to $\mc{H}_2$ is negligible. In fact, we show that our methodology is a trade-off between $\mc{H}_2$ and $\mc{H}_\infty$ designs in terms of convergence rate, where our observer converges faster than $\mc{H}_2$, and asymptotic error, where our observer provides smaller error than $\mc{H}_\infty$.

\subsection{Organization of the paper}
The rest of the paper is organized as follows. The problem is formulated in Section~\ref{sec_prob}. The main algorithm to solve the formulated problem is presented in Section~\ref{sec_mainalgo}. We also demonstrate the computational limitations of the main algorithm when dealing with large-scale network systems. Thus, Section~\ref{sec_relax} provides a structural relaxation of average state observer design to achieve computational tractability. A sufficient condition on the stabilizability of average state observer is also established in this section. Then, Section~\ref{sec_modifiedalgo} presents a modified algorithm under structural relaxation and Section~\ref{sec_sim} presents simulation results. Finally, Section~\ref{sec_conclusion} provides concluding remarks.

\section{Problem Formulation} \label{sec_prob}

%

\subsection{Clustered network system}

Consider a network represented by a digraph $\mc{G}=(\mc{V},\mc{E})$ with the set of nodes $\mc{V}$ and the set of edges ${\mc{E}\subseteq \mc{V}\times \mc{V}}$, where $(i,j)\in\mc{E}$ is an edge directed from node~$j$ to $i$ and, at time $t\in\bb{R}_{\geq 0}$, the state of each node $i\in\mc{V}$ is denoted by $x_i(t)\in\bb{R}$. The nodes are of two types: measured nodes $\mc{V}_1=\{\mu_1,\dots,\mu_m\}$, whose states are respectively measured by $m$ dedicated sensors, and unmeasured nodes $\mc{V}_2=\{\nu_1,\dots,\nu_n\}$, whose states are not measured. 
Without loss of generality, we suppose $\mc{I}_{\mc{V}_1} = \{1,\dots,m\}$ and $\mc{I}_{\mc{V}_2} = \{m+1,\dots,m+n\}$ to be the index sets of measured and unmeasured nodes. Moreover, $m\ll n$ due to limited number of sensors.

The unmeasured nodes $\mc{V}_2$ are partitioned into $k$ clusters $\mc{C}_1,\dots,\mc{C}_k$ such that, for $\alpha=1,\dots,k$ and $\beta\in\{1,\dots,k\}\setminus\{\alpha\}$,
$\mc{C}_\alpha\cap\mc{C}_\beta = \emptyset$ and $\cup_{\alpha=1}^k \mc{C}_\alpha = \mc{V}_2$.
Let $\mc{Q}=\{\mc{C}_1,\dots,\mc{C}_k\}$ denote the clustering (or partition) of $\mc{V}_2$ and $
\mf{C}_{n,k} = \{X\in\{0,1\}^{n\times k}: X\mb{1}_k = \mb{1}_n\}
$ the set of characteristic matrices of all clusterings with $k$ clusters of $n$ nodes, where a characteristic matrix $Q\in\mf{C}_{n,k}$ is defined as $[Q]_{i\alpha}=1$ if $\nu_i\in\mc{C}_\alpha$ and $[Q]_{i\alpha}=0$ otherwise.
The matrix $Q^+ = (Q^\T Q)^{-1} Q^\T$ is the left pseudo-inverse of $Q$, i.e., $Q^+ Q = I_k$, and is given by
\be \label{eq:clust_mat_Q+}
[Q^+]_{\alpha i} = \left\{\ba{cl}
\disp \frac{1}{n_\alpha}, & \text{if}~\nu_i \in \mc{C}_\alpha \\ [0.75em]
0, & \text{otherwise}
\ea\right.
\ee
where $n_\alpha = |\mc{C}_\alpha|$ is the number of nodes in $\mc{C}_\alpha$.

The clustered network system over $\mc{G}$ with measured nodes $\mc{V}_1$ and a clustering $\mc{Q}$ of unmeasured nodes $\mc{V}_2$ is defined as
\[
\bs{\Sigma}_{\mc{V}_1,\mc{Q}} : \left\{\ba{ccl}
\dot{\mb{x}}(t) &=& A\mb{x}(t) + B \mb{u}(t) \\
\mb{y}(t) &=& C\mb{x}(t)
\ea\right.
\]
where $\mb{u}(t)\in\bb{R}^p$ is the input, $\mb{x}(t) = [\colsep=2pt\ba{cc} \mb{x}_1^\T(t) & \mb{x}_2^\T(t) \ea]^\T\in\bb{R}^{m+n}$ is the state with $\mb{x}_1(t)=[\colsep=2pt\ba{ccc} x_1(t) & \dots & x_m(t) \ea]^\T$ the state of measured nodes and $\mb{x}_2(t)=[\colsep=2pt\ba{ccc} x_{m+1}(t) & \dots & x_{m+n}(t) \ea]^\T$ the state of unmeasured nodes, and $\mb{y}(t)=\mb{x}_1(t)\in\bb{R}^m$ is the output.
The state matrix $A\in\bb{R}^{(m+n)\times (m+n)}$ is Metzler, namely
\be \label{metzler_A}
\left\{\ba{ll} 
[A]_{ij}>0, & \text{if}~(i,j)\in\mc{E} \\ [0.2em]
[A]_{ii}\leq 0, & \text{if}~i=j \\ [0.2em]
[A]_{ij}=0, & \text{otherwise}.
\ea\right.
\ee
Corresponding to the partition of nodes into measured $\mc{V}_1$ and unmeasured $\mc{V}_2$, we have the following block structure of system matrices
\[
\ba{ll}
A=\left[\ba{cc}
A_{11} & A_{12}\\
A_{21} & A_{22}\ea\right], & B = \left[\ba{c}
B_1 \\ B_2
\ea\right] \\ [1em]
C=\left[\ba{cc}
I_m & 0_{m\times n}
\ea\right] &
\ea
\]
where $A_{11}\in\bb{R}^{m\times m}$, $A_{22}\in\bb{R}^{n\times n}$, $A_{12}\in\bb{R}_{\geq 0}^{m\times n}$, $A_{21}\in\bb{R}_{\geq 0}^{m\times n}$, $B_1\in\bb{R}^{m\times p}$, and $B_2\in\bb{R}^{n\times p}$.

\begin{Assumption} \label{assump1}
We assume $\rank(A_{12}) = m$.
\end{Assumption}

This assumption is reasonable because the sensors are usually placed strategically to maximize the coverage of a network system.
To interpret this assumption, suppose $A_{12}$ to be a structured matrix with a fixed zero-pattern and arbitrary non-zero values. Then, for the full-row rank of $A_{12}$, it is sufficient that (i)~no measured node is disconnected from the unmeasured nodes and (ii)~no pair of measured nodes has the same set of unmeasured nodes as their in-neighbors. 
Violating (i) means that the row of $A_{12}$ corresponding to the disconnected measured node is zero.
On the other hand, violating (ii) means that there exists, although of Lebesgue measure 0, a set of non-zero elements of the corresponding rows $i,j$ of $A_{12}$ such that the rows are linearly dependent.

\subsection{Projected network system}

By aggregating the clusters in $\mc{Q}$, one projects the state of $\bs{\Sigma}_{\mc{V}_1,\mc{Q}}$ on a lower-dimensional state space and obtains a projected network system, which provides the dynamics of the average states of clusters.
In other words, let
\[
z_\alpha(t) = \frac{1}{n_\alpha} \sum_{i\in\mc{C}_\alpha} x_i(t)
\]
be the average (mean) state of cluster $\mc{C}_\alpha$, for $\alpha=1,\dots,k$, then the average state vector $\mb{z}_\a(t) = [\colsep=2pt\ba{ccc} z_1(t) & \dots & z_k(t) \ea]^\T\in\bb{R}^k$ of $\mc{Q}$ is given by
\[
\mb{z}_\a(t) = Q^+ \mb{x}_2(t)
\]
where $Q^+$ is defined in \eqref{eq:clust_mat_Q+}.
Let $\mb{z}(t) = [\colsep=2pt\ba{cc} \mb{x}_1^\T(t) & \mb{z}_\a^\T(t) \ea]^\T$ be the projected state vector, then the projected network system can be represented as
\[
\mr{\bs{\Sigma}}_{\mc{V}_1,\mc{Q}}:\left\{\ba{ccl}
\dot{\mb{z}}(t) &=& E\mb{z}(t) + F\bs{\sigma}(t) + G\mb{u}(t) \\
\mb{0}_k &=& Q^+ \bs{\sigma}(t) \\
\mb{y}(t) &=& H\mb{z}(t)
\ea\right.
\]
where
\be \label{eq:sigma}
\bs{\sigma}(t) = (I_n - QQ^+)\mb{x}_2(t)
\ee
is the average deviation vector whose $i$-th entry, for $\nu_i\in\mc{C}_\alpha$ and $\alpha\in\{1,\dots,k\}$, is given by $[\bs{\sigma}(t)]_i = [\mb{x}_2(t)]_i - z_\alpha(t)$. The system matrices of $\mr{\bs{\Sigma}}_{\mc{V}_1,\mc{Q}}$ have the following block structure
\[
\ba{cclcl}
E &=& \left[\ba{cc}
E_{11} & E_{12} \\ E_{21} & E_{22}
\ea\right] &=& \left[\ba{cc}
A_{11} & A_{12}Q \\
Q^+  A_{21} & Q^+  A_{22}Q
\ea\right]
\\ [1em]
F &=& \left[\ba{c} F_1 \\ F_2 \ea\right] &=& \left[\ba{c}
A_{12} \\ Q^+  A_{22}
\ea\right] \\ [1em]
G &=& \left[\ba{c} G_1 \\ G_2 \ea\right] &=& \left[\ba{c} B_1 \\ Q^+  B_2 \ea\right] \\ [1em]
H &=& \left[\ba{cc} H_1 & H_2 \ea\right] &=& \left[\ba{cc} I_m & 0_{m\times k} \ea\right].
\ea
\]

\subsection{Average state observer}

The average state observer is a system that utilizes the model of projected network system to estimate the average states of clusters. Following \cite{niazi2020tcns}, the average state observer is given as
\[
\bs{\Omega}_{\mc{V}_1,\mc{Q}}:\left\{\ba{ccl}
\dot{\mb{w}}(t) &=& M_L \mb{w}(t) + K_L \mb{y}(t) + N_L \mb{u}(t) \\
\hat{\mb{z}}_\a(t) &=& \mb{w}(t) + L \mb{y}(t)
\ea\right.
\]
where
\be \label{eq:observer_matrices}
\ba{lcl}
M_L &=& E_{22} - L E_{12} \\
K_L &=& E_{21} - L E_{11} + M_L L \\
N_L &=& G_2 - L G_1
\ea
\ee
and $L\in\bb{R}^{k\times m}$ is a matrix to be designed. 

Let the estimation error be 
$
\bs{\zeta}(t) = \mb{z}_\a(t) - \hat{\mb{z}}_\a(t)
$
then
\be \label{eq:zeta}
\dot{\bs{\zeta}}(t) = M_L \bs{\zeta}(t) + R_L \bs{\sigma}(t)
\ee
where 
\be \label{eq:R_L}
R_L = F_2 - L F_1
\ee
and $M_L = R_L Q$. 

\subsection{Problem statement}

Consider the transfer matrix
\[
\mb{T}_{L,Q}(s) = (sI_k - M_L)^{-1} R_L
\]
from $\bs{\sigma}$ to $\bs{\zeta}$ in \eqref{eq:zeta}. Then, find $L\in\bb{R}^{k\times m}$ and $Q\in\mf{C}_{n,k}$ such that
\be\label{prob}
\left.\ba{cl}
\disp \min_{L,Q} & \| \mb{T}_{L,Q} \|_2 \\
\text{subject to} & M_L ~\text{is Hurwitz}
\ea\right\}
\ee
where
$
\|\mb{T}_{L,Q}\|_2 = \sqrt{\trace(W_L)}
$
is the $\mc{H}_2$ norm with
\[
W_L = \int_0^\infty \exp(M_L t) R_L R_L^\T \exp(M_L^\T t) dt
\]
the controllability gramian of $(M_L,R_L)$ in \eqref{eq:zeta}. 

\section{Main Algorithm} \label{sec_mainalgo}

In this section, we present the main algorithm for solving the problem defined in the previous section. The minimization problem \eqref{prob} has two decision variables $L$ and $Q$ with the cost 
\[
\mc{J}(L,Q) := \|\mb{T}_{L,Q}\|_2.
\]
Moreover, \eqref{prob} is convex in the decision variable $L$ but it is non-convex in $Q$ because of being mixed-integer type. That is, for a fixed $Q$, one is able to obtain the optimal design $L$ of average state observer through the LMI approach \cite{boyd2004}; however, for a fixed $L$, finding the optimal clustering $Q$ is NP-hard, and hence it is only feasible to achieve a local minimum \cite{burer2012}.

Thereby, the main algorithm is designed based on a cyclic coordinate descent scheme, which is summarized below:
\begin{enumerate}
	\item Initialization: To initialize a clustering $\mc{Q}_0$ of $n$ unmeasured nodes with $k$ clusters, generate $q=[\colsep=1pt\ba{ccc} q_1 & \dots & q_n \ea]$ with $1\leq q_i\leq k$ being a random integer such that, for every $j\in\{1,\dots,k\}$, there exists $i\in\{1,\dots,n\}$ satisfying $q_i=j$. Then, assign ${[Q_0]_{i,q_i}\leftarrow 1}$ for $i=1,\dots,n$.
	\item Repeat
	\begin{enumerate}[(i)]
		\item Let $Q=Q_0$ and find the optimal $L^*$ from \eqref{prob_LMI}.
		\item Let $L=L^*$ and compute the cost $\mc{J}(L,Q)$.
		\item Find a suboptimal $Q^*$ by Algorithm~\ref{algo:Q} and compute the cost $\mc{J}(L,Q^*)$.
		\item If $\mc{J}(L,Q^*) < \mc{J}(L,Q)$, then let $Q_0=Q^*$ and continue the loop; otherwise, return $L$ and $Q$ and stop the loop.
	\end{enumerate}
	Until convergence or maximum number of iterations. \hfill $\lrcorner$
\end{enumerate}

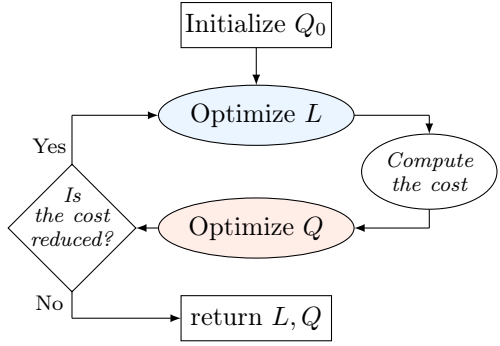
\begin{figure}[!]
\begin{center}
\begin{tikzpicture}
\draw (-0.8,-0.6) rectangle (1.2,0);
\node at (0.2,-0.3) {Initialize $Q_0$};
\draw[fill=blue2!10] (0.2,-1.5) ellipse (1.3cm and 0.4cm);
\node at (0.2,-1.5) {Optimize $L$};
\draw[fill=redd2!10] (0.2,-3) ellipse (1.3cm and 0.4cm);
\node at (0.2,-3) {Optimize $Q$};
\draw (2.5,-2.25) ellipse (0.9cm and 0.5cm);
\node[text width = 1.5cm, align = center] at (2.5,-2.25) {\scriptsize\it Compute \\ \vspace{-2pt} the cost};
\draw[rotate around = {45:(-2.25,-3)}] (-2.85,-3.6) rectangle (-1.65,-2.4);
\node[text width = 1.5cm, align = center] at (-2.25,-2.85) {\scriptsize\it Is \\ \vspace{-2pt} the cost \\ \vspace{-2pt} reduced?};
\draw (-0.8,-3.9) rectangle (1.2,-4.5);
\node at (0.2,-4.2) {return $L,Q$};

\node at (-2.55,-1.9) {\scriptsize Yes};
\node at (-2.55,-4) {\scriptsize No};

\draw (2.5,-2.75) -- (2.5,-3);
\draw (1.5,-1.5) -- (2.5,-1.5);
\draw (-2.25,-2.17) -- (-2.25,-1.5);
\draw (-2.25,-3.83) -- (-2.25,-4.2);

\path[-latex]
	(0.2,-0.6) edge (0.2,-1.1)
	(2.5,-1.5) edge (2.5,-1.75)
	(2.5,-3) edge (1.5,-3)
	(-1.1,-3) edge (-1.4,-3)
	(-2.25,-1.5) edge (-1.1,-1.5)
	(-2.25,-4.2) edge (-0.8,-4.2)
;
\end{tikzpicture}
\end{center}
\caption{The block scheme of the main algorithm.}
\label{figMainAlgo}
\end{figure}

The block scheme of the algorithm is illustrated in Figure~\ref{figMainAlgo}. In the following two subsections, we provide more details of the above algorithm.

\subsection{Algorithm to find optimal design matrix $L$}

Given $Q$, the optimal design $L^*$ of the average state observer can be obtained by solving the following LMI problem, see \cite[Chapter 9]{duan2013}, 
\begingroup
\fontsize{9.2pt}{10pt}
\be \label{prob_LMI}
\left.\colsep=1pt\ba{cl}
\min & \qquad \rho \\ [0.5em]
\text{subject to} & \left\{\ba{l}
\left[\colsep=3pt\ba{cc} \sym(X E_{22} - W E_{12}) & X \bar{F}_2 - W \bar{F}_1 \\ (X \bar{F}_2 - W \bar{F}_1)^\T & -I_n \ea\right]<0 \\ [1em]
\left[\ba{cc} -P & I_k \\ I_k & -X \ea\right]<0 \\ [1em]
\trace(P)<\rho
\ea\right.
\ea\right\}
\ee
\endgroup
where $X=X^\T\in\bb{R}^{k\times k}$, $W\in\bb{R}^{k\times m}$, and $P=P^\T\in\bb{R}^{k\times k}$ are the decision variables of the LMI problem, and $\bar{F}_i=F_i(I_n-QQ^+)$ for $i\in\{1,2\}$. The optimal design matrix is then given by
\[
L^* = X^{-1} W
\]
which minimizes the $\mc{H}_2$ norm $\|\mb{T}_{L,Q}\|_2 = \sqrt{\rho}$.

\subsection{Algorithm to find suboptimal clustering matrix $Q$}

Given $L$, a suboptimal clustering $Q^*$ can be obtained by Algorithm~\ref{algo:Q}.
At every iteration of a while loop, the algorithm consecutively assigns each unmeasured node $\nu_i\in\mc{V}_2$ to a cluster yielding the minimum cost. The while loop stops until convergence up to a prescribed tolerance level or if the maximum number of iterations are reached.

Since the algorithm depends on the initial characteristic matrix $Q_0$ and the clustering problem is non-convex, mixed-integer type, it may converge to a local minimum. To improve the results, it is recommended to repetitively run Algorithm~\ref{algo:Q} with different $Q_0$ and choose the result yielding the least cost.

\begin{algorithm}[!]
	\caption{Suboptimal clustering}\label{algo:Q}
	\begin{algorithmic}[1]
		\REQUIRE Initial characteristic matrix $Q_0$, design matrix $L$, and other matrices required to compute $\mc{J}(Q):=\mc{J}(L,Q)$
		\ENSURE Suboptimal characteristic matrix $Q^*$ 
		
		\STATE Compute the initial cost $\psi_0 = \mc{J}(Q_0)$, let $\psi_{\min} \leftarrow \psi_0$
		
		\REPEAT
		
		\FOR{$i = 1, 2, \dots, n$}
		
		\STATE Let $\beta$ be such that $\nu_i \in \mathcal{C}_\beta$, let $\psi \leftarrow \psi_{\min}$, $\alpha \leftarrow \beta$
		
		\IF{$|\mc{C}_\beta|>1$}
			\FOR{$\theta = 1 : k$, $\theta \ne \beta$}
		  		\STATE Move node~$\nu_i$ from $\mathcal{C}_\beta$ into $\mathcal{C}_\theta$
		  		
		  		\STATE Compute $Q$ and the cost $\psi_\theta = \mc{J}(Q)$
	  		
	  			\IF{$\psi_\theta < \psi$}
	  				\STATE Assign $\psi \leftarrow \psi_\theta$ and  $\alpha \leftarrow \theta$
	  			\ENDIF
	  		
	  			\STATE Move node $i$ back to the cluster $\mathcal{C}_\beta$
			\ENDFOR
		
			\STATE Move node $i$ from $\mathcal{C}_\beta$ to $\mathcal{C}_\alpha$
			
			\STATE Assign $\psi_{\min} \leftarrow \psi$ and compute $Q$
		
		\ENDIF		
		\ENDFOR
	 
	 	\STATE Assign $Q^*\leftarrow Q$
	 	
		\UNTIL convergence or the maximum number of iterations
		
		\RETURN $Q^*$.
	\end{algorithmic}
\end{algorithm}

\subsection{Computational limitation for large-scale systems} \label{subsec_complexity}

The number of nodes, especially the unmeasured nodes $n$, in large-scale network systems can be very large, precisely on the order of $10^3$ or larger. In such a case, the main algorithm will require a huge storage memory, and solving the LMI problem \eqref{prob_LMI} iteratively to obtain the optimal design matrix for each suboptimal clustering becomes computationally intractable. For example, when applying the SEDUMI solver to solve \eqref{prob_LMI}, the computational complexity is $\mc{O}(\eta^2 \theta^{2.5} + \theta^{3.5})$, \cite{labit2002}, where $\eta = 2k^2 +mk$ is the number of variables and $\theta=4k+n$ is the total number of rows in all the LMIs. This means that if $n$ is of the order $10^3$, the computational complexity is of the order $10^{9.5}$. 
Therefore, solving \eqref{prob_LMI} iteratively in the while loop of the main algorithm (Figure~\ref{figMainAlgo}) could not be practical.

\begin{figure}[!tbh]
    \centering
    \includegraphics[width=0.475\textwidth]{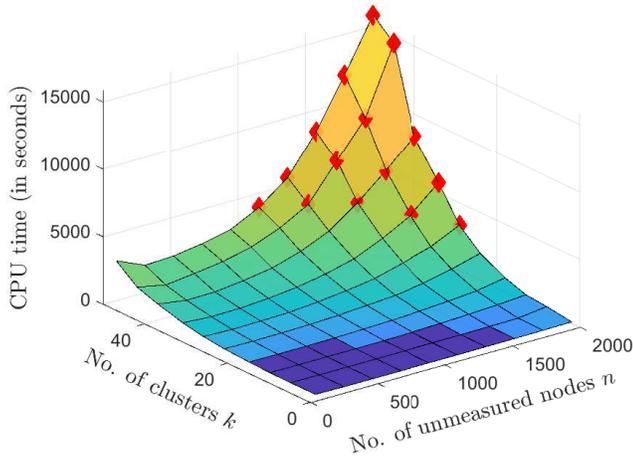}
    \caption{CPU time for solving the LMI problem \eqref{prob_LMI} as a function of the number of unmeasured nodes and clusters. The red markers show the points where the computation time exceeds an hour (i.e., 3600 seconds).}
    \label{figLMIComplexity}
\end{figure}

Figure~\ref{figLMIComplexity} shows the time it takes to solve \eqref{prob_LMI} once using SEDUMI solver in MATLAB R2021a with processor Intel Core i7 $\sim$ 3.00GHz. In this experiment, we generate a random graph $\mc{G}$ of $m+n$ nodes and consider the state matrix $A=-\mc{L}(\mc{G})$, where $\mc{L}(\mc{G})$ is the Laplacian matrix of $\mc{G}$ and the number of measured nodes $m = \lceil n / 10 \rfloor$. The characteristic matrix $Q\in\mf{C}_{n,k}$ of a clustering is chosen randomly using the method in the first step of the main algorithm. Then, after computing the matrices $E_{12},E_{22},F_1,F_2$, we solve \eqref{prob_LMI} for $n\in[100,2000]$ and $k\in[2,50]$. Notice that for $n> 1000$ and $k>25$, the computation time exceeds an hour. Particularly, when $n=2000$ and $k=50$, the computation time for solving \eqref{prob_LMI} exceeds $4.4$ hours ($15,828$ seconds). This means that if the maximum number of iterations in the main algorithm are $\eta_{\text{main}}$, then the time it takes to find optimal $L$ and $Q$ using the main algorithm is greater than $4.4 \eta_{\text{main}}$ hours plus the product of the computation time of Algorithm~\ref{algo:Q} and $\eta_{\text{main}}$. Therefore, the efficiency of the main algorithm based on the LMI computation and Algorithm~\ref{algo:Q} is not attractive for handling large-scale networks. Can we obtain a practical solution in a much more efficient manner? This problem motivates our work in the subsequent sections. 

\section{Structural Relaxation of the Design Matrix for Computational Tractability} \label{sec_relax}

As discussed in the previous section, obtaining the optimal design matrix $L$ in the main algorithm by solving the LMI problem \eqref{prob_LMI} is computationally intractable for large-scale network systems. 
In this section, we therefore propose a structural relaxation in the design matrix $L$, which achieves computational tractability and yields a suboptimal solution. 
This structural relaxation can be used in the iterations of the main algorithm to find suboptimal design matrix $L$ of the average state observer for a given characteristic matrix $Q$ of a suboptimal clustering.


\subsection{Tunability of the average state observer}
The average state observer $\bs{\Omega}_{\mc{V}_1,\mc{Q}}$ is said to be \textit{tunable} if it estimates the average state $\mb{z}_\a(t)$ at any specified rate \cite{niazi2020tcns}. Precisely, for any $\gamma>0$ and $\bs{\zeta}(0)=\bs{\zeta}_0\in\bb{R}^k$ with $\|\bs{\zeta}_0\|\leq r < \infty$, there exist a design matrix $L:=L_\gamma$ and an increasing positive-valued function $a(r)$ such that the estimation error $\bs{\zeta}(t)$ satisfies $\|\bs{\zeta}(t)\| \leq a(r) e^{-\gamma t}$.
Let
\be \label{design_L}
L = (Q^+ A_{22} - V Q^+) A_{12}^\-
\ee
where $V\in\bb{R}^{k\times k}$ is an arbitrary Hurwitz matrix.

\begin{Theorem}[see \cite{niazi2020tcns}]
\label{thm:AOiff}
The average state observer $\bs{\Omega}_{\mc{V}_1,\mc{Q}}$ is tunable \textit{if and only if}
\be \label{eq:AOiff}
\rank\left(\left[\ba{c} A_{12} \\ Q^+ A_{22} \\ Q^+ \ea\right]\right) = \rank(A_{12}).
\ee
Moreover, if \eqref{eq:AOiff} is satisfied, then, for any Hurwitz $V$, the design matrix $L$ given in \eqref{design_L} ensures that the estimation error satisfies $\dot{\bs{\zeta}}(t)= V \bs{\zeta}(t)$.
\end{Theorem}
%

The condition \eqref{eq:AOiff} applies on the structure of clustered network system $\bs{\Sigma}_{\mc{V}_1,\mc{Q}}$ with clustering $\mc{Q}$. To provide a graph-theoretic interpretation, notice that
\[
\rank\left(\left[\ba{c} A_{12} \\ Q^+ \ea\right]\right) = \rank(A_{12})
\]
is necessary for \eqref{eq:AOiff}. Since none of the columns of $Q^+$ is zero, therefore, to satisfy the above rank condition, it is necessary that none of the columns of $A_{12}$ are zero either, which means that each unmeasured node has at least one measured node as an out-neighbor. Although it can be satisfied for specific cases such as scale-free networks when their hubs are measured \cite{niazi2020tcns}, this condition is quite restrictive for real-world applications. Therefore, instead of achieving $\lim_{t\rightarrow\infty} \|\bs{\zeta}(t)\| = 0$ at an arbitrary exponential rate, we aim to minimize $\limsup_{t\rightarrow\infty} \|\bs{\zeta}(t)\|$. However, as mentioned earlier, achieving the optimal solution for this problem is computationally expensive in the main algorithm (Figure~\ref{figMainAlgo}). Thus, we aim for a suboptimal solution by exploiting the structure of the average deviation vector.

\subsection{Structural relaxation of the design matrix}
The average deviation vector $\bs{\sigma}(t)$ acts as a structured unknown input in both the projected network system $\mathring{\bs{\Sigma}}_{\mc{V}_1,\mc{Q}}$ and the estimation error equation \eqref{eq:zeta}. It is structured because $\bs{\sigma}(t)\in\Ker(Q^+)$ and it is unknown because it is a function of the unmeasured state $\mb{x}_2(t)$ defined in \eqref{eq:sigma}. 

Theorem~\ref{thm:AOiff} implies that if the average state observer $\bs{\Omega}_{\mc{V}_1,\mc{Q}}$ is tunable, then, for any Hurwitz $V\in\bb{R}^{k\times k}$, choosing the design matrix $L$ as in \eqref{design_L} ensures that $R_L\bs{\sigma}\equiv \mb{0}_k$ and the estimation error $\bs{\zeta}(t)$ exponentially converges to zero as $t\rightarrow\infty$ at an arbitrary rate $\gamma>0$. If $\bs{\Omega}_{\mc{V}_1,\mc{Q}}$ is not tunable, then $R_L\bs{\sigma}\equiv \mb{0}_k$ if the clustered network system $\bs{\Sigma}_{\mc{V}_1,\mc{Q}}$ is average detectable \cite[Theorem 5]{niazi2020tcns}, which requires regularity in the inter-cluster and intra-cluster topologies of $\bs{\Sigma}_{\mc{V}_1,\mc{Q}}$. In the case of intra-cluster consensus or synchronization in $\bs{\Sigma}_{\mc{V}_1,\mc{Q}}$, the signal $R_L\bs{\sigma}(t)$ converges to zero asymptotically, ensuring the convergence of estimation error \eqref{eq:zeta} to zero if $L$ is such that $M_L = R_L Q$ is Hurwitz (see \cite[Theorem 7]{niazi2020tcns}).

Generally speaking, the conditions ensuring $R_L\bs{\sigma}\equiv \mb{0}_k$ or $R_L\bs{\sigma}(t)\rightarrow \mb{0}_k$ as $t\rightarrow\infty$ are quite restrictive. Therefore, it is reasonable to find $L$ that minimizes $\|R_L\bs{\sigma}(t)\|$ by exploiting the structure of $\bs{\sigma}(t)$. Notice that, for a given $V\in\bb{R}^{k\times k}$, the design matrix $L$ given by \eqref{design_L} is the least-square solution to $R_L = VQ^+$ minimizing
$
\min_{L\in\bb{R}^{k\times m}} \|R_L - VQ^+\|
$
which, by \eqref{eq:equivalence}, is equivalent to
$
\min_{L\in\bb{R}^{k\times m}} \|R_L\bs{\sigma}(t)\|.
$
Thus, fixing $L$ as in \eqref{design_L} gives $R_L$ as a function of $V$, i.e.,
\[
R_V := R_{L} = Q^+ A_{22}(I_n - A_{12}^\- A_{12}) + V Q^+ A_{12}^\- A_{12}.
\]
Then, we find optimal $V$ as follows.

\begin{Lemma} \label{lemma:V_star}
Consider $L$ as in \eqref{design_L}.
Then,
\be \label{eq:V*}
V^* = Q^+ A_{22} Q
\ee
is the minimizing solution to
$
\min_{V\in\bb{R}^{k\times k}} \|R_V - VQ^+\|.
$
\end{Lemma}
\begin{pf}
See Appendix~\ref{appendix:prooflemma3}. \qed
\end{pf}

The choice of $L$ as in \eqref{design_L} with $V=V^*=Q^+ A_{22} Q$ minimizes the effect of average deviation in \eqref{eq:zeta}. However, such a choice of $V$ may not ensure the stability of average state observer, which is characterized by $M_L = R_L Q$ being Hurwitz. In the next subsection, we will show that by perturbing $V$ by a scalar parameter $\phi\in\bb{R}$, we can stabilize the average state observer under some mild sufficient conditions.

\subsection{Stabilizability of the average state observer under structural relaxation of the design matrix} \label{subsec_stab}


Consider $\mc{G}_\nu = (\mc{V}_2,\mc{E}_\nu)$ to be the induced subgraph formed by the unmeasured nodes $\mc{V}_2$, where $\mc{E}_\nu = \mc{E}\cap(\mc{V}_2\times\mc{V}_2)$. The off-diagonal entries of the matrix $A_{22}$ constitute the edge configuration of $\mc{G}_\nu$. Therefore, the subgraph $\mc{G}_\nu$ is weakly connected if and only if the matrix $A_{22}+A_{22}^\T$ is irreducible. Note that a matrix is said to be reducible if it can be transformed to a block upper-triangular form by simultaneous row/column permutations---otherwise, it is said to be irreducible.

The weak connectivity of the induced subgraph $\mc{G}_\nu$ can also be established by considering its undirected version $\ol{\mc{G}}_\nu$, where the edges of $\ol{\mc{G}}_\nu$ are obtained by ignoring the directions from the edges of $\mc{G}_\nu$. Then, $\mc{G}_\nu$ is weakly connected if and only if $\ol{\mc{G}}_\nu$ is connected, which is equivalent to having the rank of its Laplacian matrix $\rank(\mc{L}(\ol{\mc{G}}_\nu)) = n-1$. 

For simplicity of notation, let $[A_{22}]_{ij} = a_{ij}$, for $i,j=1,\dots,n$. Then, for every unmeasured node $\nu_i\in\mc{V}_2$,
\be \label{eq:s_i}
s_i = \sum_{\substack{j=1 \\ j\neq i}}^n \left( a_{ij} + a_{ji} \right)
\ee
is the sum of the weights of all edges going into and emerging from $i$ within the induced subgraph $\mc{G}_\nu$. That is, the sum of the weights of all edges of the node~$\nu_i$'s in-neighbors and out-neighbors. Finally, recall from \eqref{metzler_A} that all the diagonal entries of $A_{22}$ are non-positive, i.e., $[A_{22}]_{jj} = a_{jj}\leq 0$ for $j=1,\dots,n$.

\begin{Assumption} \label{assump2}
We assume the following:
\begin{enumerate}[(i)]
\item The induced subgraph $\mc{G}_\nu=(\mc{V}_2,\mc{E}_\nu)$ is weakly connected.
\item For every unmeasured node $\nu_i\in\mc{V}_2$, $s_i \leq 2 |a_{ii}|$, and there exists at least one $\nu_j\in\mc{V}_2$ such that $s_{j} < 2 |a_{jj}|$.
\end{enumerate}
\end{Assumption}

In what follows, we provide a sufficient condition for the stabilizability of the average state observer $\bs{\Omega}_{\mc{V}_1,\mc{Q}}$ under the structural relaxation of design matrix $L$ as in \eqref{design_L}, where 
\be \label{V_phi}
V =: V_\phi =\phi Q^+ A_{22} Q
\ee
is perturbed by a scalar $\phi\in\bb{R}$. Similarly, for a fixed $Q\in\mf{C}_{n,k}$, we use the notations 
$R_\phi := R_L$ and $M_\phi := M_L$.
The stability of the average state observer is achieved by showing that, for some range of $\phi$, the state matrix $M_\phi = R_\phi Q$ is stabilized, or made Hurwitz. Precisely, we say that the average state observer is \textit{stabilizable} by the gain parameter $\phi\in\bb{R}$ if, for every $Q\in\mf{C}_{n,k}$, there exists $\psi\in\bb{R}$ such that, for every $\phi>\psi$, the matrix $M_\phi=R_\phi Q$ is Hurwitz.

\begin{Theorem} \label{thm:Mr_stability}
Let Assumption~\ref{assump1} and \ref{assump2} hold. Then, the average state observer $\bs{\Omega}_{\mc{V}_1,\mc{Q}}$ of a clustered network system $\bs{\Sigma}_{\mc{V}_1,\mc{Q}}$ is stabilizable by the gain parameter $\phi\in\bb{R}$ if the characteristic matrix $Q\in\mf{C}_{n,k}$ is such that $\rank(A_{12} Q) = k$, where $k$ is the number of clusters.
\end{Theorem}

The proof of this theorem is provided in Appendix~\ref{appendix:proofthm4}. Here, we briefly discuss the implications of this theorem for the clustering problem.

\begin{Corollary} \label{coro:k_leq_m}
If the characteristic matrix $Q\in\mf{C}_{n,k}$ is such that $\rank(A_{12}Q)=k$, then the number of clusters~$k$ must be less than or equal to the number of measured nodes~$m$.
\end{Corollary}
\begin{pf}
Assume the contrary that $\rank(A_{12} Q)=k$ and ${k>m}$. We know $\rank(A_{12})=m$ and $\rank(Q)=k$, and
\[\ba{ccl}
\rank(A_{12}Q)&\leq&\min(\rank(A_{12}),\rank(Q))\\
&=&\min(m,k).
\ea\]
Thus, $\rank(A_{12}Q)\leq m < k$, which is a contradiction. \qed
\end{pf}

This means that if we employ Theorem~\ref{thm:Mr_stability} in the clustering algorithm to ensure stabilizability, then the number of clusters of unmeasured nodes cannot exceed the number of measured nodes.

Denote the neighbor set of measured nodes $\mc{V}_1$ with respect to the unmeasured nodes $\mc{V}_2$ as $\mc{N}_{\mc{V}_1\leftarrow\mc{V}_2}\subseteq\mc{V}_2$, contains all the unmeasured nodes that have at least one measured out-neighbor, i.e.
\be \label{eq:NV1V2}
\mc{N}_{\mc{V}_1\leftarrow\mc{V}_2} = \{\nu_j\in\mc{V}_2: \exists\,(\mu_i,\nu_j)\in\mc{E} ~\text{for some}~ \mu_i\in\mc{V}_1\}.
\ee

\begin{Corollary} \label{coro:rank(A_12 Q)}
If the characteristic matrix $Q\in\mf{C}_{n,k}$ is such that $\rank(A_{12}Q)=k$, then, for every $\alpha\in\{1,\dots,k\}$, we have $\mc{C}_\alpha \cap \mc{N}_{\mc{V}_1\leftarrow\mc{V}_2} \neq \emptyset$.
\end{Corollary}
\begin{pf}
Assume the contrary that $\rank(A_{12}Q)=k$ and there exists some $\alpha\in\{1,\dots,k\}$ such that $\mc{C}_\alpha\cap\mc{N}_{\mc{V}_1\leftarrow\mc{V}_2} = \emptyset$.
Let $A_{12} = \left[\colsep=2pt\ba{ccc}
\mb{a}_1 & \dots & \mb{a}_n
\ea\right]$, where $\mb{a}_1,\dots,\mb{a}_n\in\bb{R}^m$ are the columns of $A_{12}\in\bb{R}_{\geq 0}^{m\times n}$. Then, we can write
$
A_{12}Q=\left[\begin{array}{ccc}
\mb{p}_1 & \dots & \mb{p}_k
\end{array}\right]
$
where $\mb{p}_1,\dots,\mb{p}_k\in\bb{R}^{m}$ are the columns of $A_{12}Q$ with 
$
\mb{p}_\alpha = \sum_{j\in\mc{C}_\alpha} \mb{a}_j
$
for $\alpha=1,\dots,k$.
Since $\rank(A_{12}Q)=k$, we have that, for every $\alpha\in\{1,\dots,k\}$, $\mb{p}_\alpha\neq \mb{0}_{m}$.
On the other hand, since we assumed that there exists $\alpha\in\{1,\dots,k\}$ such that ${\mc{C}_\alpha\cap\mc{N}_{\mc{V}_1\leftarrow\mc{V}_2} = \emptyset}$, which means the cluster $\mc{C}_\alpha$ does not contain any node from the neighbor set $\mc{N}_{\mc{V}_1\leftarrow\mc{V}_2}$. That is, $\mb{a}_j=0$ for all $j\in\mc{C}_\alpha$, therefore $\mb{p}_\alpha = 0$. This implies that $\rank(A_{12}Q)<k$, which is a contradiction. \qed
\end{pf}

Therefore, to ensure that the sufficient condition of stabilizability is satisfied, it is necessary to choose clusters such that $k\leq m$ and that every cluster contains at least one unmeasured node in the neighbor set of the measured nodes.

\subsection{Sufficiency of the clustering constraint for stabilizability in a generic sense}

The clustering constraint of Corollary~\ref{coro:rank(A_12 Q)} 
\[
\mc{C}_\alpha \cap \mc{N}_{\mc{V}_1\leftarrow\mc{V}_2} \neq \emptyset, \; \forall \alpha\in\{1,\dots,k\}
\]
is necessary for satisfying the condition $\rank(A_{12} Q)=k$ of Theorem~\ref{thm:Mr_stability} for given matrices $A_{12}$ and $Q$. However, instead of the rank, if we consider a generic rank \cite{lin1974, murota1987}, where only the non-zero pattern of $A_{12}$ is taken into account, then this condition can also be proven to be sufficient.

The generic rank of a matrix, denoted as $\grank(\cdot)$, is defined as the maximum rank among all choices of the non-zero entries of the matrix. 
In general, the rank of a matrix is less than or equal to its generic rank. However, for any matrix, whose non-zero entries are chosen randomly in some interval, its rank is equal to the generic rank almost always (i.e., with probability 1), except for the non-zero entries of the matrix in some proper algebraic variety, which is of Lebesgue measure zero \cite{dion2003}.

We can represent the non-zero pattern of any $Z\in\bb{R}^{m\times n}$ by a bipartite graph $\mc{G}_Z=(\mc{V}_r,\mc{V}_c,\mc{E}_Z)$, where $\mc{V}_r=\{r_1,\dots,r_m\}$ is the index set of the rows of $Z$, $\mc{V}_c = \{c_1,\dots,c_n\}$ is the index set of the columns of $Z$, and $\mc{E}_Z\subseteq\mc{V}_r\times \mc{V}_c$ is the set of edges defined as $(r_i,c_j)\in\mc{E}_Z$ if $[Z]_{r_i c_j}\neq 0$. A {\em matching} in a bipartite graph $\mc{G}_X$ is the set of edges such that no two edges have a vertex in common, whereas a {\em maximum matching} is a matching with the maximum possible number of edges \cite{godsil2001}.
Then, from \cite{liu2011}, we know that \textit{the generic rank of $Z\in\bb{R}^{m\times n}$ is equal to the size of maximum matching in $\mc{G}_Z$}.

\begin{Theorem} \label{thm:grank}
Let $Q\in\mf{C}_{n,k}$ be the characteristic matrix of some clustering $\mc{Q}$ of $n$ unmeasured nodes $\mc{V}_2$ and let Assumption~\ref{assump1} hold. Further, assume $k\leq m$, where $m$ is the number of measured nodes. Then, $\grank(A_{12}Q)=k$ {\em if and only if} the clustering $\mc{Q}=\{\mc{C}_1,\dots,\mc{C}_k\}$ is such that, for every $\alpha\in\{1,\dots,k\}$, $\mc{C}_\alpha \cap \mc{N}_{\mc{V}_1\leftarrow\mc{V}_2} \neq \emptyset$. 
\end{Theorem}
\begin{pf}[Proof of sufficiency]
Assume that the clustering $\mc{Q}=\{\mc{C}_1,\dots,\mc{C}_k\}$ is such that, for every $\alpha\in\{1,\dots,k\}$, it holds $\mc{C}_\alpha\cap\mc{N}_{\mc{V}_1\leftarrow\mc{V}_2}\neq \emptyset$. By Assumption~\ref{assump1}, we have $\rank(A_{12})=m$, and, since $\rank(A_{12})\leq \grank(A_{12})\leq m$, $\grank(A_{12}) = m$, which implies that a maximum matching of the bipartite graph $\mc{G}_{A_{12}}$ is of size $m$. That is, for all the $m$ measured nodes there are distinct $m$ neighbors in $\mc{N}_{\mc{V}_1\leftarrow\mc{V}_2}$, which implies $|\mc{N}_{\mc{V}_1\leftarrow\mc{V}_2}|\geq m$. Since $k\leq m$ and $|\mc{N}_{\mc{V}_1\leftarrow\mc{V}_2}|\geq m$, we have $|\mc{N}_{\mc{V}_1\leftarrow\mc{V}_2}|\geq k$. Let $\nu_1,\dots,\nu_k\in\mc{N}_{\mc{V}_1\leftarrow\mc{V}_2}$ to be $k$ unmeasured nodes that are in the clusters $\mc{C}_1,\dots,\mc{C}_k$, respectively. Then, a matching of size $k$ of the bipartite graph $\mc{G}_{A_{12}}$ is
\begin{center}
\begin{tikzpicture}[scale=0.6]
\vertex[white, fill=black, minimum width = 4.5mm] (mu1) at (0,3) {\scriptsize $\mu_1$};
\vertex[white, fill=black, minimum width = 4.5mm] (mu2) at (3,3) {\scriptsize $\mu_2$};
\node at (5,3) {$\dots$};
\vertex[white, fill=black, minimum width = 4.5mm] (muk) at (7,3) {\scriptsize $\mu_k$};
\node at (9,3) {$\dots$};
\vertex[white, fill=black, minimum width = 4.5mm] (mum) at (11,3) {\scriptsize $\mu_m$};

\vertex[fill=red!40!, minimum width = 3mm] (nu1) at (1,0) {\scriptsize $\nu_1$};
\vertex[fill=red!40!, minimum width = 2.5mm] at (1.7,0) {};
\node at (2.25,0) {\scriptsize $\dots$};
\vertex[fill=red!40!, minimum width = 2.5mm] at (2.8,0) {};
\draw[dashed, black!80!] (1.9,0) ellipse (1.25cm and 0.65cm);
\node at (1.9,-1) {$\mc{C}_1$};

\vertex[fill=red!40!, minimum width = 3mm] (nu2) at (4,0) {\scriptsize $\nu_2$};
\vertex[fill=red!40!, minimum width = 2.5mm] at (4.7,0) {};
\node at (5.25,0) {\scriptsize $\dots$};
\vertex[fill=red!40!, minimum width = 2.5mm] at (5.8,0) {};
\draw[dashed, black!80!] (4.9,0) ellipse (1.25cm and 0.65cm);
\node at (4.9,-1) {$\mc{C}_2$};

\node at (6.75,0) {$\dots$};

\vertex[fill=red!40!, minimum width = 3mm] (nuk) at (8,0) {\scriptsize $\nu_k$};
\vertex[fill=red!40!, minimum width = 2.5mm] at (8.7,0) {};
\node at (9.25,0) {\scriptsize $\dots$};
\vertex[fill=red!40!, minimum width = 2.5mm] at (9.8,0) {};
\draw[dashed, black!80!] (8.9,0) ellipse (1.25cm and 0.65cm);
\node at (8.9,-1) {$\mc{C}_k$};

\path[very thick]
	(mu1) edge (nu1)
	(mu2) edge (nu2)
	(muk) edge (nuk)
;

\end{tikzpicture}
\end{center}
where $\mu_1,\dots,\mu_m$ are the $m$ measured nodes representing the rows of $A_{12}$ and the nodes on the right side are the $n$ unmeasured nodes representing the columns of $A_{12}$. The unmeasured nodes are partitioned into clusters $\mc{C}_1,\dots,\mc{C}_k$. The bipartite graph $\mc{G}_{A_{12}Q}$ is obtained by aggregating the clusters in $\mc{G}_{A_{12}}$. Then, from the matching of size $k$ illustrated above for $\mc{G}_{A_{12}}$, we obtain a maximum matching of size $k$ for $\mc{G}_{A_{12}Q}$ 
\begin{center}
\begin{tikzpicture}[scale=0.6]
\vertex[white, fill=black, minimum width = 4.5mm] (mu1) at (0,3) {\scriptsize $\mu_1$};
\vertex[white, fill=black, minimum width = 4.5mm] (mu2) at (3,3) {\scriptsize $\mu_2$};
\node at (5,3) {$\dots$};
\vertex[white, fill=black, minimum width = 4.5mm] (muk) at (7,3) {\scriptsize $\mu_k$};
\node at (9,3) {$\dots$};
\vertex[white, fill=black, minimum width = 4.5mm] (mum) at (11,3) {\scriptsize $\mu_m$};

\vertex[fill=red!40!, minimum width = 4.5mm] (nu1) at (1,0) {\scriptsize $c_1$};

\vertex[fill=red!40!, minimum width = 4.5mm] (nu2) at (4,0) {\scriptsize $c_2$};

\node at (6,0) {$\dots$};

\vertex[fill=red!40!, minimum width = 4.5mm] (nuk) at (8,0) {\scriptsize $c_k$};

\path[very thick]
	(mu1) edge (nu1)
	(mu2) edge (nu2)
	(muk) edge (nuk)
;

\end{tikzpicture}
\end{center}
where the clusters $\mc{C}_1,\dots,\mc{C}_k$ are represented as super nodes $c_1,\dots,c_k$, respectively. Thus, we have $\grank(A_{12}Q)=k$.

{\em Proof of necessity.}
See the proof of Corollary~\ref{coro:rank(A_12 Q)}. \qed
\end{pf}


In Section~\ref{subsec_stab}, stabilizability of the average state observer $\bs{\Omega}_{\mc{V}_1,\mc{Q}}$ was defined as the existence of some $\psi\in\bb{R}$ such that $M_\phi = R_\phi Q$ is Hurwitz for all $\phi>\psi$.
Theorem~\ref{thm:Mr_stability} provides a sufficient condition for the stabilizability of average state observer as the full-column rank of $A_{12}Q$. Similarly, the average state observer is stabilizable in a \textit{generic sense} if $A_{12}Q$ has full column generic rank. This means that there exists $\psi\in\bb{R}$ {\em almost always} such that $M_{\phi,Q}$ is Hurwitz for every $\phi>\psi$. The term `almost always' indicates that, for any submatrix $A_{12}\in\bb{R}_{\geq 0}^{m\times n}$ belonging to a clustered network system $\bs{\Sigma}_{\mc{V}_1,\mc{Q}}$, the rank $A_{12}Q$ is equal to $k$ with probability one if the condition of Theorem~\ref{thm:grank} is satisfied.

\begin{Corollary} \label{thm:H2OptClust}
Let Assumption~\ref{assump1} and \ref{assump2} hold. Then, for any clustering $\mc{Q}$ of unmeasured nodes with $k\leq m$ clusters, the average state observer $\bs{\Omega}_{\mc{V}_1,\mc{Q}}$ is stabilizable in a generic sense if $\forall\alpha\in\{1,\dots,k\}$, $\mc{C}_\alpha\cap\mc{N}_{\mc{V}_1\leftarrow\mc{V}} \neq \emptyset$.
\end{Corollary}

The proof follows directly from Theorem~\ref{thm:Mr_stability} and \ref{thm:grank}.

The sufficient condition in Corollary~\ref{thm:H2OptClust} corresponds to the clustering $\mc{Q}$, which needs to be satisfied by the clustering algorithm for ensuring stabilizability in a generic sense.

\section{Modified Algorithm under Structural Relaxation of the Design Matrix} \label{sec_modifiedalgo}

Under the structural relaxation of the design matrix $L$ as in \eqref{design_L} with $V$ given in \eqref{V_phi}, the problem~\eqref{prob} is modified as follows
\begingroup
\fontsize{10pt}{10pt}
\be \label{prob:modified}
\left.\colsep=1pt\ba{cl}
\disp \min_{\rho>0, Q\in\mf{C}_{n,k}} & \disp \mc{J}(\phi,Q) := \trace(W_{\phi,Q}) \\ [1em]
\text{subject to} & \left\{\ba{l}
M_{\phi,Q} ~\text{is Hurwitz} \\
\mc{C}_\alpha \cap \mc{N}_{\mc{V}_1\leftarrow\mc{V}_2} \neq \emptyset, ~ \forall\alpha\in\{1,\dots,k\}
\ea\right.
\ea\right\}
\ee
\endgroup
where $\mc{N}_{\mc{V}_1\leftarrow\mc{V}_2}$ is the neighbor set of measured nodes defined in \eqref{eq:NV1V2} and
\[
W_{\phi,Q} = \int_0^\infty \exp(M_{\phi,Q} t) R_{\phi,Q} R_{\phi,Q}^\T \exp(M_{\phi,Q}^\T t) dt
\]
with $M_{\phi,Q} = R_{\phi,Q} Q$ and 
\[
R_{\phi,Q} = Q^+ A_{22} \left( I_n - (I_n - \phi Q Q^+ ) A_{12}^\- A_{12} \right).
\]
Note that $M_{\phi,Q},R_{\phi,Q}$ are the matrices $M_L,R_L$ in \eqref{eq:zeta}.

%
%
%

The modified algorithm is summarized below:
\begin{enumerate}
    \item Initialization: Initialize a clustering $\mc{Q}$ with the characteristic matrix $Q_0\in\mf{C}_{n,k}$ using lines~\ref{algo:clust_rho_p11}--\ref{algo:clust_rho_p12} of  Algorithm~\ref{algo:clust_rho}.
    \item Repeat
        \begin{enumerate}[(i)]
            \item Let $Q=Q_0$ and find optimal $\phi^*$ using Algorithm~\ref{algo:rho_star}.
            \item Let $\phi=\phi^*$ and compute the cost $\mc{J}(\phi,Q)$.
            \item Find a suboptimal $Q^*$ using Algorithm~\ref{algo:clust_rho} and compute $\mc{J}(\phi,Q^*)$.
            \item If $\mc{J}(\phi,Q^*) < \mc{J}(\phi,Q)$, then let $Q_0=Q^*$ and continue the loop; otherwise, return $\phi$ and $Q$, and stop the loop.
        \end{enumerate}
    Until convergence or maximum number of iterations.
\end{enumerate}

Similar to the main algorithm, the modified algorithm also follows the scheme of Figure~\ref{figMainAlgo}, where instead of $L$, we optimize the scalar gain parameter $\phi$.

\subsection{Algorithm to find optimal gain parameter $\phi$}

For a fixed $Q\in\mf{C}_{n,k}$, the cost in \eqref{prob:modified} is simply written as $\mc{J}(\phi)$ and $M_\phi:=M_{\phi,Q}$.
Then, the problem of finding the optimal gain parameter is defined as follows: Find $\phi^*\in\bb{R}$ such that
\be \label{eq:H2Optproblem}
\phi^* = \arg\min_{\phi\in\bb{R}} \mc{J}(\phi) \quad \text{subject to}~M_\phi~\text{is Hurwitz}.
\ee
Note that the problem \eqref{eq:H2Optproblem} is a convex optimization problem with a single decision variable, see e.g., \cite[Chapter 3]{boyd1994} and \cite[Chapter 4]{boyd2004}. Therefore, a global minimum can be achieved easily by a simple algorithm as Algorithm~\ref{algo:rho_star}.

	
\begin{algorithm}[!]
\caption{Incremental search algorithm}
\label{algo:rho_star}
\begin{algorithmic}[1]
	\REQUIRE Matrices required to compute $M_{\phi}$ and $\mc{J}(\phi)$, tolerance $\ul{\varepsilon}>0$, parameter $\eta\geq 2$, and initial step size $\varepsilon>0$
	\ENSURE Optimal solution $\phi^*$ to Problem~\eqref{eq:H2Optproblem}
	
	\STATE Initialize $\phi<0$ such that $M_\phi$ is not Hurwitz
	
	\STATE Assign $\varepsilon_1 \leftarrow \varepsilon$
	
	\REPEAT
		\STATE Compute $M_\phi$
		\IF{$M_\phi$ is Hurwitz}
			\STATE Assign $\psi\leftarrow\phi$, $\phi\leftarrow\phi-\varepsilon_1$, $\varepsilon_1\leftarrow\varepsilon_1/\eta$
		\ELSE
			\STATE Assign $\phi\leftarrow\phi + \varepsilon_1$
		\ENDIF
	\UNTIL $\eta\varepsilon_1 \leq \ul{\varepsilon}$
	
	\STATE Assign $\phi\leftarrow \psi + \varepsilon_1 \eta$ and compute the cost $c=\mc{J}(\phi)$
	
	\REPEAT
		\STATE Assign $\phi\leftarrow\phi+\varepsilon$ and compute  $c_1 = \mc{J}(\phi)$
		\IF{$c_1 > c$}
			\STATE Assign $\phi\leftarrow\phi - 2\varepsilon$ and $\varepsilon\leftarrow\varepsilon/\eta$
			
			\STATE Compute $c=\mc{J}(\phi)$
		\ELSE
			\STATE Assign $c\leftarrow c_1$
		\ENDIF

	\UNTIL $\eta\varepsilon \leq \ul{\varepsilon}$
	\RETURN $\phi^*\leftarrow \phi + \varepsilon \eta$.
\end{algorithmic}
\end{algorithm}

The main idea of the above algorithm is to initialize $\phi\in\bb{R}$ and continue to increment it with a small $\varepsilon>0$ in order to search for the optimal solution. The value of $\varepsilon>0$ is initialized arbitrarily and then, in the algorithm, is reduced iteratively by dividing it with parameter $\eta\geq 2$. This reduction achieves the required tolerance level $\ul{\varepsilon}>0$ towards the actual optimal solution $\phi^*$. In the algorithm, whenever $\phi$ passes the optimal value, we define a smaller interval around that optimal value, divide the interval into several points, choose $\varepsilon$ to be the length of these divisions, and search for the optimal solution in this interval. This process is done iteratively until a required tolerance level is achieved.

In the first part of Algorithm~\ref{algo:rho_star}, we find the minimum $\psi>0$ such that, for $\phi=\psi$, we have $M_\phi$ Hurwitz. Then, in the second part, we initialize $\phi=\psi+\varepsilon_1\eta$, where $\eta \varepsilon_1\leq \ul{\varepsilon}$ is the achieved tolerance level, and increment it by $\varepsilon>0$ until we pass the optimal solution, which is the global minimum. This is because before the global minimum was reached, the cost $\mc{J}(\phi)$ in non-increasing at every iteration. However, when the cost increases at a certain iteration, it indicates that $\phi$ has surpassed the global minimum. At this point, we know that the solution lies in the interval $[\phi-2\varepsilon,\phi]$. Therefore, we decrement $\phi$ by $2\varepsilon$, decrease the value of $\varepsilon$ by dividing it by $\eta$, and restart the search process in the specified interval. This process is repeated until the solution $\phi^*$ is within the specified tolerance $\ul{\varepsilon}$ to the true optimal value.

\begin{figure}[!]
    \centering
    \includegraphics[width=0.475\textwidth]{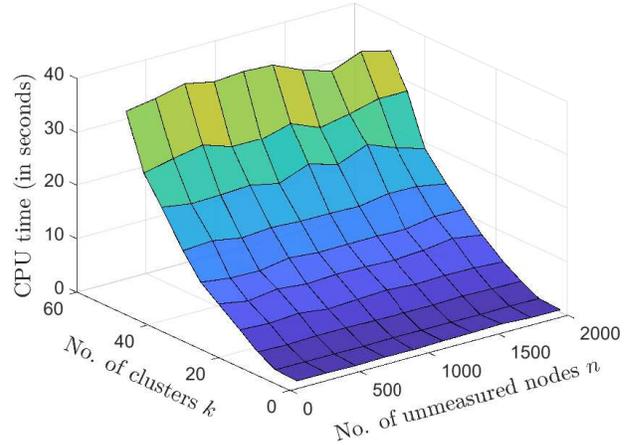}
    \caption{CPU time of Algorithm~\ref{algo:rho_star} as a function of the number of unmeasured nodes and clusters.}
    \label{figRhostarComplexity}
\end{figure}

For the initial step size $\varepsilon=10^{-3}$, reducing parameter $\eta=10$, and tolerance $\ul{\varepsilon}=10^{-9}$, Figure~\ref{figRhostarComplexity} shows the computation time of Algorithm~\ref{algo:rho_star} in MATLAB R2021a with processor Intel Core i7 $\sim$ 3.00GHz. The setup of this experiment is the same as that of Section~\ref{subsec_complexity}.
Notice that the CPU time grows almost linearly in $k$ and stays constant with respect to $n$. Particularly, for number of nodes $n=2000$ and clusters $k=50$, the computation time of Algorithm~\ref{algo:rho_star} is approximately $30$ seconds, which is quite feasible to be implemented iteratively in the modified algorithm.

\subsection{Modified algorithm to find suboptimal clustering}

For a fixed $\phi\in\bb{R}$, the problem is to find a clustering $\mc{Q}=\{\mc{C}_1,\dots,\mc{C}_k\}$ with characteristic matrix $Q\in\mf{C}_{n,k}$ such that the cost $\mathcal{J}(Q)$ in problem \eqref{prob:modified} is minimized subject to the constraint on clusters $\mc{C}_\alpha\cap\mc{N}_{\mc{V}_1\leftarrow\mc{V}_2}\neq \emptyset$, for $\alpha=1,\dots,k$. 

Algorithm~\ref{algo:clust_rho} finds a suboptimal clustering solution to \eqref{prob:modified} by using a greedy approach. In the first part of the algorithm, lines~\ref{algo:clust_rho_p11}--\ref{algo:clust_rho_p12}, we initialize the $k$ clusters such that the second constraint of \eqref{prob:modified} is satisfied.
Let $\mc{S}_1,\dots,\mc{S}_k$ be $k$ clusters of a subset of $\mc{V}_2$ and define
\[
\mc{N}_{\mc{S}_\alpha\leftrightarrow\mc{V}_2} = \{j\in\mc{V}_2:(i,j)\in\mc{E}_\nu ~\text{or}~(j,i)\in\mc{E}_\nu, ~\text{for}~i\in\mc{S}_\alpha\}
\]
to be the set of in-neighbors and out-neighbors of $\mc{S}_\alpha$.
Notice that, by Assumption~\ref{assump1}, we have $|\mc{N}_{\mc{V}_1\leftarrow\mc{V}_2}|\geq m$ and, by Corollary~\ref{coro:k_leq_m}, $m\geq k$. Therefore, one can find the non-empty cluster subsets $\mc{S}_1,\dots,\mc{S}_k$ in line~\ref{algo:clust_rho_p11}. Also, by Assumption~\ref{assump2}, the induced subgraph $\mc{G}_\nu$ is weakly connected. Therefore, the while loop in lines~\ref{algo:clust_rho_while1}--\ref{algo:clust_rho_while2} that iteratively traverses the graph $\mc{G}_\nu$ using the breadth-first search to include the immediate neighbors of each subset $\mc{S}_1,\dots,\mc{S}_k$ terminates, where the subsets $\mc{S}_1,\dots,\mc{S}_k$ are disjoint and their union is equal to the set of unmeasured nodes $\mc{V}_2$. Since each $\mc{S}_1,\dots,\mc{S}_k$ were initialized by partitioning the neighbor set $\mc{N}_{\mc{V}_1\leftarrow\mc{V}_2}$, the initial clustering $\mc{Q}_0=\{\mc{C}_1,\dots,\mc{C}_k\}$ obtained in line~\ref{algo:clust_rho_p12} satisfies the stabilizability constraint of \eqref{prob:modified}.
Finally, the second part of Algorithm~\ref{algo:clust_rho} iteratively moves all the unmeasured nodes that are not in the neighbor set to a cluster yielding the minimum cost. The algorithm stops when the specified tolerance level $\delta>0$ is reached.







\begin{algorithm}[!]
\caption{Suboptimal clustering algorithm under the stabilizability constraint on clusters}
\label{algo:clust_rho}
\begin{algorithmic}[1]
\REQUIRE Matrices needed to compute $\mc{J}(Q)$, neighbor set $\mc{N}_{\mc{V}_1\leftarrow\mc{V}_2}$, and tolerance $\delta>0$ (e.g., $10^{-6}$)
\ENSURE Suboptimal clustering $\mc{Q}=\{\mc{C}_1,\dots,\mc{C}_k\}$

\STATE Move each $j\in\mc{N}_{\mc{V}_1\leftarrow\mc{V}_2}$ to one of the subsets $\mc{S}_1,\dots,\mc{S}_k$ such that, $\forall\alpha\in\{1,\dots,k\}$, $\mc{S}_\alpha\neq\emptyset$
\label{algo:clust_rho_p11}

\REPEAT \label{algo:clust_rho_while1}
	\STATE Assign $\mc{S}_1 \leftarrow \mc{S}_1 \cup (\mc{N}_{\mc{S}_1\leftrightarrow\mc{V}_2} \setminus \mc{N}_{\mc{V}_1\leftarrow\mc{V}_2})$
	\FOR{$\alpha=2,\dots,k$}
		\STATE Assign $\mc{S}_\alpha \leftarrow \mc{S}_\alpha \cup (\mc{N}_{\mc{S}_1\leftrightarrow\mc{V}_2} \setminus \mc{N}_{\mc{V}_1\leftarrow \mc{V}_2} \setminus \mc{S}_{\alpha-1})$
	\ENDFOR

\UNTIL $\mc{S}_1\cup\dots\cup\mc{S}_k=\mc{V}_2$
\label{algo:clust_rho_while2}

\STATE Assign $\mc{C}_1\leftarrow\mc{S}_1, \mc{C}_2\leftarrow\mc{S}_2, \dots,\mc{C}_k\leftarrow\mc{S}_k$

\STATE Assign $\mc{Q}_0\leftarrow\{\mc{C}_1,\dots,\mc{C}_k\}$ and compute $c_0 = \mc{J}(Q_0)$
\label{algo:clust_rho_p12}

\STATE Assign $\mc{Q}_1\leftarrow\mc{Q}_0$

\REPEAT
	\STATE Assign $c_1\leftarrow c_0$
	\FOR{$i\in\mc{V}_2\setminus\mc{N}_{\mc{V}_1\leftarrow\mc{V}_2}$}
		\STATE Assign $\mc{Q}_2\leftarrow\mc{Q}_1$
		\STATE Let $\beta$ be such that $i\in\mc{C}_\beta$
		\FOR{$\alpha=1,\dots,k$ and $\alpha\neq\beta$}
			\STATE Move $i$ to $\mc{C}_\alpha$ and update $\mc{Q}_2$ accordingly
			\STATE Compute $c_2=\mc{J}(Q_2)$
			\IF{$c_2<c_0$}
				\STATE Assign $c_0\leftarrow c_2$ and $\mc{Q}_1\leftarrow\mc{Q}_2$
			\ELSE
				\STATE Move $i$ back to $\mc{C}_\beta$ and $\mc{Q}_2\leftarrow\mc{Q}_1$
			\ENDIF
		\ENDFOR	
	\ENDFOR
	\STATE Assign $\mc{Q}\leftarrow\mc{Q}_1$
\UNTIL $c_1-c_0<\delta$, i.e., specified tolerance to convergence
\RETURN $\mc{Q}=\{\mc{C}_1,\dots,\mc{C}_k\}$.
\end{algorithmic}
\end{algorithm}

\section{Simulation Results} \label{sec_sim}

\subsection{An example of linear flow on directed networks}
\label{subsec_flow}

Linear flow networks model many real-world large-scale infrastructures such as urban traffic networks, electrical power grids, and water distribution systems \cite{kaiser2021}. The state $x_i(t)$ of each node~$i$ represents an amount of some physical quantity at time~$t$, which evolves according to
\be \label{eq:flow}
\dot{x}_i(t) = \sum_{j\in\mc{N}_i^{\text{in}}} a_{ij} x_j(t) - \sum_{h\in\mc{N}_i^{\text{out}}} a_{ji} x_i(t) + \sum_{g=1}^p b_{ig} u_g(t)
\ee
where $\mc{N}_i^{\text{in}}$ and $\mc{N}_i^{\text{out}}$ are the sets of $i$'s in-neighbors and out-neighbors in the directed graph $\mc{G}=(\mc{V},\mc{E})$, respectively. The first term on the right hand side of the above equation represents the total inflow to $i$ from its in-neighbors, the second term represents the total outflow from $i$ to its out-neighbors, and the third term represents external inputs. In vector-form, the model is written as $
\dot{\mb{x}}(t) = A \mb{x}(t) + B\mb{u}(t)
$
where $A = \mc{A}(\mc{G}) - \mc{D}_{\text{out}}(\mc{G})$ with $\mc{D}_{\text{out}} = \diag(\mb{1}^\T \mc{A}(\mc{G}))$ the out-degree matrix and $\mc{A}(\mc{G})$ the adjacency matrix of the directed graph $\mc{G}$.

For evaluating our clustering-based average state observer design on a linear flow network, we suppose the number of unmeasured nodes $n=1000$, measured nodes $m=10$, clusters $k=10$, and inputs $p=50$. The graph $\mc{G}$ with $n+m$ nodes is generated using the Erd\H{o}s-R\'{e}nyi random graph model with the probability of a directed edge between each pair of nodes chosen to be $P_{\text{edge}}=0.05$ and the weight of the edge chosen uniformly randomly in $(0,1)$. The obtained directed graph is such that the induced subgraph $\mc{G}_\nu$ formed by unmeasured nodes is weakly connected because $P_{\text{edge}}>2\ln(n)/(n)$, \cite{erdos1959}. The input $u_g(t) = a_g \sin(w_g t + b_g)$, where $a_g,w_g,b_g$ are chosen uniformly randomly in the intervals $(-0.05,0.05), (0,0.5),(-\pi,\pi)$, respectively. The input matrix $B$ is generated by considering the probability that each input $u_g$ acts on node~$i$ equal to $0.01$. The initial condition $x_i(0)$ of the model \eqref{eq:flow} is chosen uniformly randomly in $(0,1)$.

\begin{figure}[!]
    \begin{center}
        \begin{tikzpicture}
        \node at (0,0) {\includegraphics[width=0.475\textwidth, trim=240 140 180 120, clip]{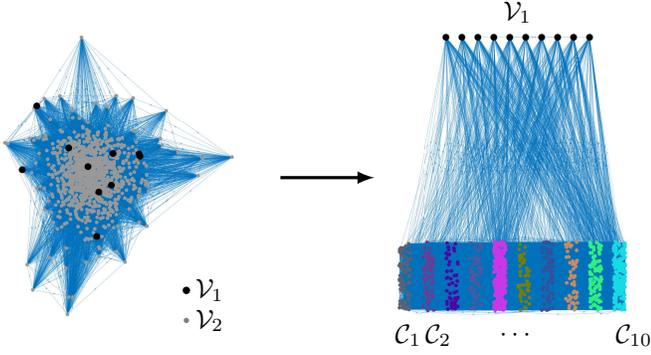}};
        \draw[-latex, very thick, black, fill=white] (-0.5,0) -- (0.75,0);
        \filldraw[color=black, fill=black] (-1.75,-1.5) circle (1.2pt);
        \filldraw[color=black!50, fill=black!50] (-1.75,-1.9) circle (0.8pt);
        \node[anchor=west] at (-1.75,-1.5) {$\mc{V}_1$};
        \node[anchor=west] at (-1.75,-1.9) {$\mc{V}_2$};
        
        \node at (2.65,2.2) {$\mc{V}_1$};
        \node at (1.2,-2.1) {$\mc{C}_1$};
        \node at (1.6,-2.1) {$\mc{C}_2$};
        \node at (4.2,-2.1) {$\mc{C}_{10}$};
        \node at (2.65,-2.1) {$\cdots$};
        \end{tikzpicture}
    \end{center}
    
    \caption{Suboptimal clustering obtained from the modified algorithm. The large black nodes depict the measured nodes, whereas the smaller colored nodes are the unmeasured nodes.}
    \label{fig:clusters}
\end{figure}

\begin{figure}[!]
\begin{center}
    \begin{tikzpicture}
    \node at (0,0) {\includegraphics[width=0.475\textwidth, trim=0 0 30 15, clip]{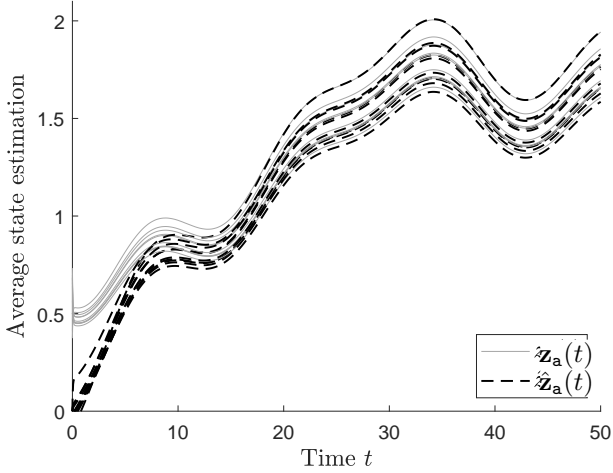}};
    \draw[white, fill=white] (3.2,-1.45) rectangle (3.85,-2.2);
    \node[anchor=west] at (3.05,-1.65) {$\mb{z}_\a(t)$};
    \node[anchor=west] at (3.05,-2.03) {$\hat{\mb{z}}_\a(t)$};
    \end{tikzpicture}
\end{center}
    \caption{Average state estimation (ASE).}
    \label{fig:avgest}
\end{figure}

\begin{figure}[!]
    \centering
    \includegraphics[width=0.475\textwidth, trim=0 0 30 15, clip]{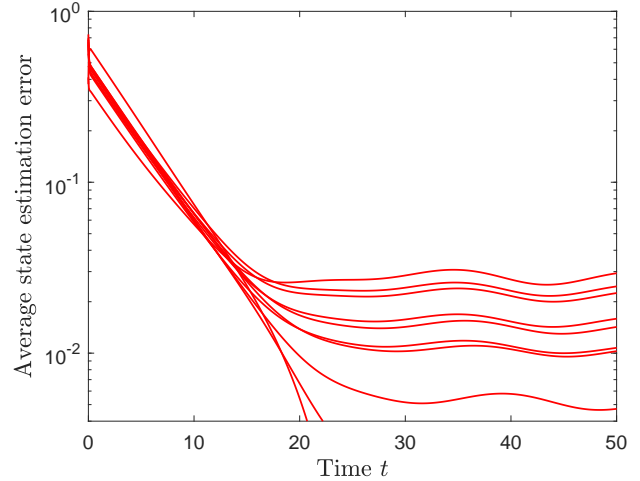}
    \caption{Evolution of the ASE error.}
    \label{fig:error}
\end{figure}

We run the modified algorithm with the maximum number of iterations equal to $10$ and tolerance equal to $10^{-8}$. In Algorithm~\ref{algo:rho_star}, the tolerance parameter $\ul{\varepsilon}=10^{-8}$, parameter $\eta=10$, and initial step size $\varepsilon=1$. In Algorithm~\ref{algo:clust_rho}, the tolerance parameter $\delta=10^{-6}$. The resulting optimal gain parameter $\phi=1.7516$ with the suboptimal clustering shown in Figure~\ref{fig:clusters}.
Then, the average state observer $\bs{\Omega}_{\mc{V}_1,\mc{Q}}$ is obtained by choosing $M_L,K_L,N_L$ as in \eqref{eq:observer_matrices}, where $L$ is in \eqref{design_L} with $V=V_\phi$ in \eqref{V_phi}. The average state estimation (ASE) result is illustrated in Figure~\ref{fig:avgest} and the ASE error in Figure~\ref{fig:error}. We obtain the percentage asymptotic estimation error
\[
\zeta_{\%} = \limsup_{t\rightarrow\infty} \frac{\|\mb{z}_\a(t) - \hat{\mb{z}}_\a(t)\|}{\|\mb{z}_\a(t)\|} \times 100
\]
to be $\zeta_{\%} \approx 0.96\%$.

\subsection{Comparison with LMI-based $\mc{H}_2$ and $\mc{H}_\infty$ observer designs}

Given the suboptimal clustering (Figure~\ref{fig:clusters}), we compare our design methodology (Algorithm~\ref{algo:rho_star}) with $\mc{H}_2$ and $\mc{H}_\infty$ designs. The $\mc{H}_2$ design is obtained by solving the LMI problem \eqref{prob_LMI}, which minimizes the $\mc{H}_2$ norm of the error system \eqref{eq:zeta}. The $\mc{H}_\infty$ design minimizes the $\mc{H}_\infty$ norm of \eqref{eq:zeta}, and is obtained by solving the following LMI problem (see \cite{duan2013}):
\begingroup
\fontsize{9.2pt}{10pt}
\[
\left.\colsep=1pt\ba{cl}
\min & \qquad \rho \\ [0.5em]
\text{subject to} & \left[\colsep=3pt\ba{ccc} \sym(X E_{22} - W E_{12}) & X \bar{F}_2 - W \bar{F}_1 & I_k \\
(X \bar{F}_2 - W \bar{F}_1)^\T & -\rho I_n & 0_{n\times k} \\
I_k & 0_{k\times n} & -\rho I_k
\ea\right]<0
\ea\right\}
\]
\endgroup
where $X=X^\T\in\bb{R}^{k\times k}$ and $W\in\bb{R}^{k\times m}$ are the decision variables. Then, similar to $\mc{H}_2$ design, the $\mc{H}_\infty$ design is given by $L=X^{-1} W$.

\begin{figure}[!]
    \centering
	\includegraphics[width=0.475\textwidth, trim=0 0 30 15, clip]{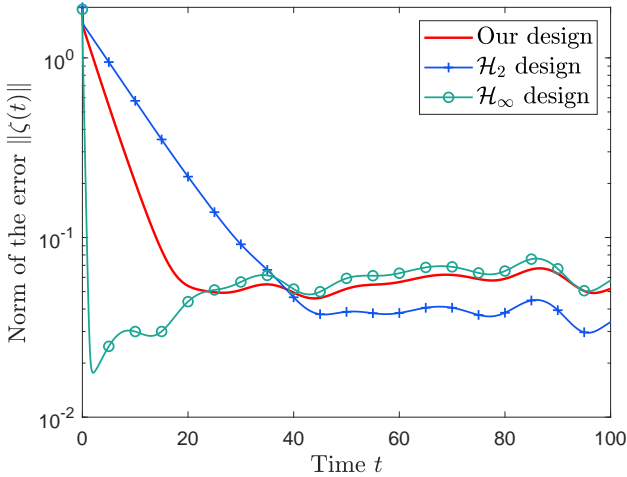}
    \caption{Norm of the ASE error for $\mc{H}_2$ and $\mc{H}_\infty$ designs, and our design methodology (Algorithm~\ref{algo:rho_star}).}
    \label{fig:H2Hinf}
\end{figure}

The norm of the average state estimation error for all three designs is illustrated in Figure~\ref{fig:H2Hinf}. As the $\mc{H}_2$ design minimizes the asymptotic estimation error, it yields the smallest error at the steady state ($\zeta_{\%}\approx 0.64\%$). On the other hand, the $\mc{H}_\infty$ design minimizes the maximum of the estimation error, which is $\bs{\zeta}(0)$ at the initial time $t=0$, therefore it yields a fastest convergence but largest error at the steady state ($\zeta_{\%}\approx 1.1\%$). Our design method provides a trade-off between the convergence rate and steady state error. The convergence times of our method is approximately $20$ seconds, which is better than $\mc{H}_2$ ($>40$ seconds) but worse than $\mc{H}_\infty$ ($<3$ seconds). However, in terms of computation time, our design gave an optimal gain parameter in less $5$ seconds, whereas $\mc{H}_2$ and $\mc{H}_\infty$ took approximately 370 and 485 seconds, respectively, to provide solutions, which is $74$ times slower than our method.

\subsection{Comparison of estimation errors for undirected Erd\H{o}s-R\'{e}nyi and Scale-free networks}

\begin{figure}[!]
    \centering
    \includegraphics[width=0.475\textwidth, trim=40 10 30 10, clip]{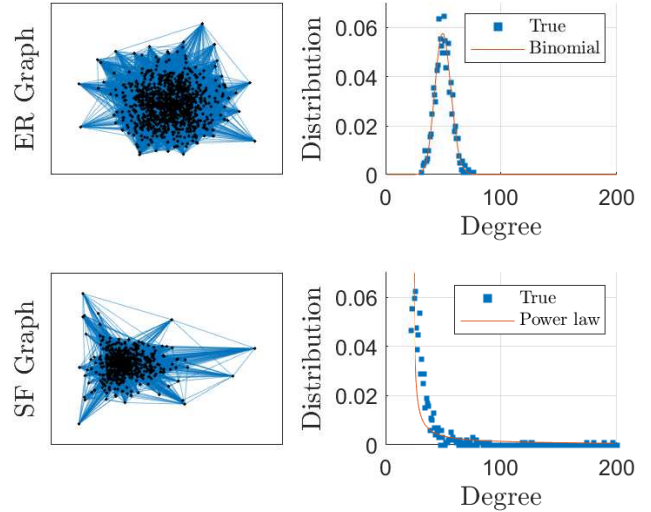}
    \caption{Erd\H{o}s-R\'{e}nyi and Scale-free graphs.}
    \label{fig:ERSFgraphs}
\end{figure}

In this experiment, we compare the effect of network degree distribution on the ASE error. We first generate 100 undirected Erd\H{o}s-R\'{e}nyi (ER) graphs with $n+m$ nodes ($n=100$, $m=5$) and edge probability $P_{\text{edge}}$ chosen uniformly random in $(0.1,0.25)$. Then, we generate 100 undirected Scale-free (SF) graphs with the same number of nodes, and bias and number of edges chosen uniformly randomly in $(2,2.5)$ and $(500,1000)$, respectively. The degree distributions of an example of these graphs is shown in Figure~\ref{fig:ERSFgraphs}. The degree distribution of an ER graph resembles a binomial distribution and that of an SF graph a power law distribution, where the accuracy increases as the number of nodes increase.

\begin{figure}[!]
    \centering
    \includegraphics[width=0.475\textwidth, trim=0 0 30 15, clip]{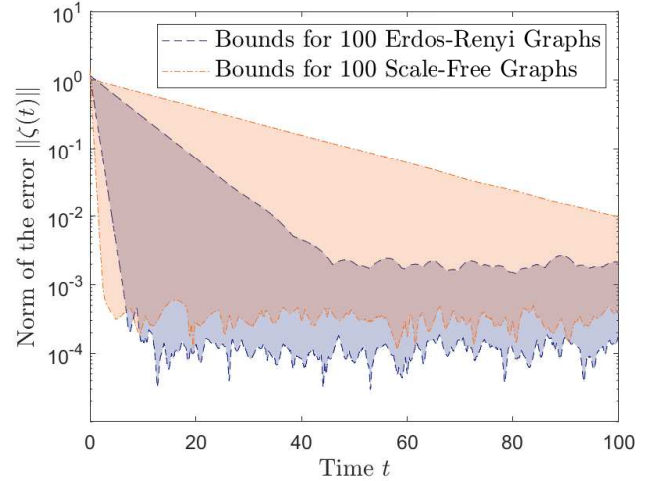}
    \caption{Bounds on the norm of the ASE error for randomly generated Erd\H{o}s-R\'{e}nyi and Scale-free graphs.}
    \label{fig:ERSFbounds}
\end{figure}

We consider a linear flow network over these graphs with input defined similarly as in Section~\ref{subsec_flow}. The number of clusters $k=5$. Figure~\ref{fig:ERSFbounds} shows the regions containing the norm of the ASE error $\|\bs{\zeta}(t)\|$ for randomly generated ER and SF graphs. First, notice that the convergence in the best scenario of an SF graph is more than twice as fast as the best scenario of an ER graph. This is because SF graphs are better suited for average observability when the hubs are taken as measured nodes \cite{niazi2020tcns}. On the other hand, the convergence in the worst scenario of an ER graph is much faster than the worst scenario of an SF graph. Also, the percentage asymptotic error $\zeta_{\%}$ in the worst scenario of ER graph is about $0.21\%$, whereas in the worst scenario of SF graph is about $1\%$, which is around four times larger. This might be due to the intra-connectivity of clusters, which is much lower in SF graphs than ER graphs. The research in this direction is however an interesting prospect.

\section{Concluding Remarks} \label{sec_conclusion}

Monitoring large-scale network systems becomes challenging when the computational and sensing resources are limited. By taking these limitations into account, we studied the problem of clustering-based average state observer design to enable aggregated monitoring of network systems through the estimation of the average states of clusters. The proposed algorithm finds a clustering and an average state observer design yielding a minimal asymptotic estimation error. As the clustering is a non-convex, mixed integer type optimization problem, we presented a greedy algorithm to obtain a suboptimal solution. On the other hand, due to the computational infeasibility of finding optimal average state observer design at every iteration, we sought a structural relaxation of its design matrix to achieve computational tractability. Such a relaxation is realizable because of the structural property of the average deviation vector, which acts as an unknown input in the dynamics of the estimation error. The design matrix can then be fixed with a single perturbation parameter that is tuned to minimize the asymptotic estimation error. Under the structural relaxation, we provided a sufficient condition for the stabilizability of average state observer, which is then incorporated in the clustering algorithm.

Although the compromise on optimality is marginal, the proposed design methodology gains a significant advantage over $\mc{H}_2$ and $\mc{H}_\infty$ average state observer designs in terms of computation time. Moreover, the transient response of our average state observer design is much faster than the $\mc{H}_2$ design, which is at the expense of larger asymptotic estimation error. Nonetheless, our design yields a smaller asymptotic estimation error than the $\mc{H}_\infty$ design, which has the fastest transient response.

Our future work includes the estimation of other aggregated state profiles such as variance and higher moments of clusters' state trajectories. Combining a sensor placement algorithm with a clustering-based average state observer design is also an interesting prospect.

\appendix



\section{Proof of Lemma~\ref{lemma:V_star}} \label{appendix:prooflemma3}
We have
$
R_V - VQ^+ = (Q^+ A_{22} - VQ^+ ) (I_n - A_{12}^\- A_{12}).
$
If $R_V-VQ^+ = 0_{k\times n}$, then
$
\Ker(Q^+ A_{22} - VQ^+) \supseteq \Ker(A_{12})
$
because the columns of $I_n - A_{12}^\- A_{12}$ form a complete basis of $\Ker(A_{12})$. This implies that 
$
Q^+ A_{22} - VQ^+ = W A_{12}
$
for some $W\in\bb{R}^{k\times m}$. However, if the ideal solution does not exist, the minimizing solution is the least-square solution 
$
V = (Q^+ A_{22} - W A_{12}) Q
$
which implies
\begingroup
\fontsize{9.5pt}{10pt}
\[
R_V - VQ^+ = (Q^+ A_{22}(I_n - QQ^+ ) + W A_{12} QQ^+ )(I_n -A_{12}^\- A_{12}).
\]
\endgroup
Finally, the minimizing solution to
\[
\min_{W\in\bb{R}^{k\times m}} \|Q^+ A_{22}(I_n - QQ^+ ) + W A_{12} QQ^+\|
\]
is 
$
W = Q^+ A_{22} (QQ^+ - I_n) QQ^+ A_{12}^+ = 0_{k\times m}
$
because $(QQ^+ - I_n) Q = 0_{n\times k}$. Thus, $V = V^* = Q^+ A_{22} Q$ is the minimizing solution to 
\[
\min_{V\in\bb{R}^{k\times k}} \|R_V - VQ^+\|.
\]

\section{Proof of Theorem~\ref{thm:Mr_stability}} \label{appendix:proofthm4}

Assume $\rank(A_{12}Q) = k$ and that Assumption~\ref{assump1} holds, i.e., $\rank(A_{12})=m$. Then,
\[\ba{ccl}
\rank(A_{12} Q) &=& \rank((A_{12} A_{12}^\-)^{-\frac{1}{2}} A_{12} Q) \\ [0.5em]
&=& \rank(Q^+ A_{12}^\- A_{12} Q) \\ [0.5em]
&=& k
\ea\]
where we used the properties $\rank(X^\T X) = \rank(X)$ and $\rank(YX) = \rank(X)$, for some matrix $X\in\bb{R}^{a\times b}$ and a non-singular $Y\in\bb{R}^{a\times a}$. This implies that the matrix
\be \label{eq:mat_S}
\ba{ccl}
S &:=& Q^\T A_{12}^\- A_{12} Q \\
&=& [(A_{12} A_{12}^\-)^{-\frac{1}{2}} A_{12} Q]^\T (A_{12} A_{12}^\-)^{-\frac{1}{2}} A_{12} Q
\ea
\ee
is positive definite because $(A_{12} A_{12}^\-)^{-\frac{1}{2}} A_{12} Q$ has full column rank. 
Recall $R_\phi:=R_L$, $L_\phi:=L$, and $V_\phi$ from \eqref{eq:R_L}, \eqref{design_L}, and \eqref{V_phi}, respectively, then we can write $M_\phi = R_\phi Q = \phi X S + Y$
where
\be \label{eq:mat_XY}
\ba{ccl}
X &=& Q^+ A_{22} Q (Q^\T Q)^{-1} \\ [0.5em]
Y &=& Q^+ A_{22}( I_n - A_{12}^\- A_{12})Q.
\ea
\ee

\begin{Lemma}[S-stability \cite{arrow1958, ostrowski1962, carlson1968}] \label{lemma:S-stability}
Let $X,S\in\bb{R}^{n\times n}$ be two matrices. If $X+X^\T$ is negative definite and $S=S^\T$ is positive definite, then the product $XS$ is Hurwitz.
\end{Lemma}

\begin{Lemma} \label{lemma:QXQ}
Let $X\in\bb{R}^{n\times n}$. If $X+X^\T$ is negative definite, then, for every non-negative matrix $Q\in\bb{R}^{n\times k}$ with $\rank(Q)=k$, the matrix $Q^\T X Q$ is Hurwitz.
\end{Lemma}
\begin{pf}
Since $X+X^\T < 0$, therefore $Q^\T X Q + Q^\T X^\T Q <0$ for a non-negative $Q\in\bb{R}^{n\times k}$ with $\rank(Q)=k$. Then, the result follows by Lyapunov's theorem (see \cite[Theorem 2.2.1]{horn1991}). \qed
\end{pf}

\begin{Lemma} \label{thm:V*stability}
If Assumption~\ref{assump2} holds, then, for every $Q\in\mf{C}_{n,k}$, the matrix $Q^+ A_{22} Q$ is Hurwitz.
\end{Lemma}
\begin{pf}
First, if Assumption~\ref{assump2}(i) holds, then the symmetric part of $A_{22}$, $\mc{S}(A_{22}) = A_{22} + A_{22}^\T$, is irreducible. That is, an undirected graph $\ol{\mc{G}}_\nu$ capturing the structure of $\mc{S}(A_{22})$ is connected. Thus, the Laplacian matrix of $\ol{\mc{G}}_\nu$ defined as
\[
[\mc{L}(\ol{\mc{G}}_\nu)]_{ij} = \left\{\ba{cl}
s_i, & \text{if}~i=j \\
-(a_{ij}+a_{ji}), & \text{if}~i\neq j
\ea\right.
\]
is of rank $n-1$ and nullity $1$, where $s_i$ is defined in \eqref{eq:s_i}. Since $\mc{L}(\ol{\mc{G}}_\nu)$ is positive semi-definite, we have, for every $\mb{v}\in\bb{R}^n$, $\mb{v}^\T \mc{L}(\ol{\mc{G}}_\nu) \mb{v} \geq 0$. Moreover, $0\in\eig(\mc{L}(\ol{\mc{G}}_\nu))$ with algebraic multiplicity $1$ because $\ol{\mc{G}}_\nu$ is connected, therefore we have $\mb{v}^\T \mc{L}(\ol{\mc{G}}_\nu) \mb{v} = 0$ if and only if $\mb{v} = a \mb{1}_n$, for $a\in\bb{R}$, i.e., in the direction of the eigenvector of $\mc{L}(\ol{\mc{G}}_\nu)$ corresponding to the $0$ eigenvalue.

Second, if Assumption~\ref{assump2}(ii) holds, then
$
\mc{S}(A_{22}) = -\mc{L}(\ol{\mc{G}}_\nu) - \mc{D}
$
where $\mc{D} = \diag(2|a_{11}| - s_1,\dots,2|a_{nn}| - s_n)$ is a diagonal matrix, which is positive semi-definite because, for all $i\in\{1,\dots,n\}$, we have $s_i \leq 2|a_{ii}|$ and, for at least one $j\in\{1,\dots,n\}$, we have $s_j < 2|a_{jj}|$. Thus, for every $\mb{v}\in\bb{R}^n$, we have $\mb{v}^\T \mc{D} \mb{v} \geq 0$. However, we know that $\mb{1}_n^\T \mc{D} \mb{1}_n > 0$ and, for some $\mb{v}_1\in\bb{R}^n$ such that $\mb{v}_1^\T \mc{D} \mb{v}_1 = 0$, we have ${\mb{v}_1^\T \mc{L}(\ol{\mc{G}}_\nu)\mb{v}_1 > 0}$. Thus, $\mc{L}(\ol{\mc{G}}_\nu)+\mc{D}$ is positive definite because $\mb{v}^\T (\mc{L}(\ol{\mc{G}}_\nu)+\mc{D}) \mb{v} > 0$ for every $\mb{v}\in\bb{R}^n$, implying that 
\be \label{eq:S(A22)}
\mc{S}(A_{22}) = A_{22} + A_{22}^\T = - (\mc{L}(\ol{\mc{G}}_\nu)+\mc{D})
\ee
is negative definite. Therefore, for any $Q\in\mf{C}_{n,k}$, we have $Q^+ A_{22} Q$ Hurwitz by Lemma~\ref{lemma:S-stability} and \ref{lemma:QXQ}. \qed
\end{pf}

Thus, by Lemma~\ref{thm:V*stability} and \ref{lemma:S-stability}, the matrix $X = Q^+ A_{22}Q (Q^\T Q)^{-1}$ is Hurwitz. Moreover, from \eqref{eq:S(A22)} it holds that 
$
X + X^\T = Q^+ (A_{22} + A_{22}^\T ) Q^{+\T} < 0.
$
Therefore, again by Lemma~\ref{lemma:S-stability} and the fact that $S$ in \eqref{eq:mat_S} is positive definite, we have that $XS$ is Hurwitz. 

\begin{Lemma} \label{lemma:X+Y}
Let $Z,Y\in\bb{R}^{k\times k}$ be any matrices with $Z$ being Hurwitz. Then, there exists $\psi\in\bb{R}$ such that, for every $\phi > \psi$, the matrix $\phi Z + Y$ is Hurwitz.
\end{Lemma}
\begin{pf}
By Lyapunov's theorem (see \cite[Theorem 2.2.1]{horn1991}), it is necessary and sufficient for $Z$ to be Hurwitz that there exists a positive definite matrix $P=P^\T$ such that $PZ + Z^\T P$ is negative definite. Since $Z$ is Hurwitz, there exists $P>0$ such that $PZ + Z^\T P<0$. For such a $P$, we have
\[
P(\phi Z + Y) + (\phi Z + Y)^\T P = \phi (PZ + Z^\T P) + (PY + Y^\T P)
\]
negative definite if, and only if, there exists $\phi$ such that, for every $\mb{v}\in\bb{R}^k$,
\be \label{eq:Hurwitz_rho1}
\phi \mb{v}^\T (PZ + Z^\T P) \mb{v} < - \mb{v}^\T (PY + Y^\T P) \mb{v}.
\ee
Thus, if the above inequality holds, then $\phi Z+Y$ is Hurwitz. In the following, we show that indeed there exists $\psi$ such that \eqref{eq:Hurwitz_rho1} is satisfied for every $\phi>\psi$.

Since $PZ+Z^\T P<0$, we have $\mb{v}^\T (PZ + Z^\T P) \mb{v} < 0$ for every $\mb{v}\in\bb{R}^n$. Therefore, dividing both sides of \eqref{eq:Hurwitz_rho1} by $\mb{v}^\T (PZ + Z^\T P) \mb{v}$ changes the sign of the inequality and gives
\be \label{eq:Hurwitz_rho2}
\phi > \frac{\mb{v}^\T (PY + Y^\T P) \mb{v}}{|\mb{v}^\T (PZ + Z^\T P) \mb{v}|}.
\ee
Let
\[
\ba{cclc}
\mb{v}_1 &=& \disp \arg\min_{\mb{v}\in\bb{R}^n} & \disp \frac{\left| \mb{v}^\T (PZ + Z^\T P) \mb{v} \right|}{\mb{v}^\T\mb{v}} \\ [1em]
\mb{v}_2 &=& \disp \arg\max_{\mb{v}\in\bb{R}^n} & \disp \frac{\mb{v}^\T (PY + Y^\T P) \mb{v}}{\mb{v}^\T\mb{v}}
\ea
\]
then choosing
\[
\psi = \frac{\mb{v}_2^\T ( PY + Y^\T P) \mb{v}_2}{|\mb{v}_1^\T (PZ + Z^\T P) \mb{v}_1 |} \frac{\mb{v}_1^\T\mb{v}_1}{\mb{v}_2^\T\mb{v}_2}.
\]
implies
\[
\psi \geq \frac{\mb{v}^\T (PY + Y^\T P) \mb{v}}{|\mb{v}^\T (PZ + Z^\T P) \mb{v}|}.
\]
Therefore, every $\phi > \psi$ ensures \eqref{eq:Hurwitz_rho2}, and thus \eqref{eq:Hurwitz_rho1}. \qed
\end{pf}

Since the matrix $XS$ is Hurwitz, therefore, by Lemma~\ref{lemma:X+Y}, there exists $\psi\in\bb{R}$ such that $M_\phi = \phi XS + Y$ is Hurwitz $\forall\phi>\psi$, where $X,Y$ are defined in \eqref{eq:mat_XY} and $S$ in \eqref{eq:mat_S}. This concludes the proof of Theorem~\ref{thm:Mr_stability}.

\bibliographystyle{plain}
\bibliography{bib_journal}

\end{document}